%% file: main.tex
\documentclass[reqno]{amsart}
\usepackage {amsmath,amssymb}
\usepackage {color,graphicx}
\usepackage [latin1]{inputenc}
\usepackage {times,mathptm}
\usepackage {float}

\input {definitions.tex}

\begin {document}

\title [Tropical Algebraic Geometry]{Tropical Algebraic Geometry}
\author {Andreas Gathmann}
\address {Andreas Gathmann, Fachbereich Mathematik, Technische Universität
  Kaiserslautern, Postfach 3049, 67653 Kaiserslautern, Germany}
\email {andreas@mathematik.uni-kl.de}

\maketitle

\input {introduction.tex}

\input {tropcurve.tex}

\input {classical.tex}

\input {enum.tex}

\input {summary.tex}

\input {biblio.tex}
\end {document}

%% file: definitions.tex
\newif\ifpdf
        \ifx\pdfoutput\undefined
        \pdffalse 
        \else
        \pdfoutput=1 
        \pdftrue
        \fi

\newcommand {\CC}{{\mathbb C}}
\newcommand {\NN}{{\mathbb N}}
\newcommand {\PP}{{\mathbb P}}
\newcommand {\QQ}{{\mathbb Q}}
\newcommand {\RR}{{\mathbb R}}
\newcommand {\ZZ}{{\mathbb Z}}

\DeclareMathOperator{\ev}{ev}
\DeclareMathOperator{\val}{val}
\DeclareMathOperator{\Val}{Val}
\DeclareMathOperator{\Div}{Div}
\DeclareMathOperator{\Pic}{Pic}
\DeclareMathOperator{\Log}{Log}
\DeclareMathOperator{\area}{area}

\newcommand {\preprint}[2]{preprint \discretionary {#1/}{#2}{#1/#2}}

\theoremstyle {definition}

\theoremstyle {remark}

%% file: introduction.tex


There are many examples in algebraic geometry in which complicated geometric or
algebraic problems can be transformed into purely combinatorial problems. The
most prominent example is probably given by toric varieties --- a certain class
of varieties that can be described purely by combinatorial data, e.g.\ by
giving a convex polytope in an integral lattice. As a consequence, most
questions about these varieties can be transformed into combinatorial questions
on the defining polytope that are then hopefully easier to solve.

Tropical algebraic geometry is a recent development in the field of algebraic
geometry that tries to generalize this idea substantially. Ideally, every
construction in algebraic geometry should have a combinatorial counterpart in
tropical geometry. One may thus hope to obtain results in algebraic geometry by
looking at the tropical (i.e.\ combinatorial) picture first and then trying to
transfer the results back to the original algebro-geometric setting.

The origins of tropical geometry date back about twenty years. One of the
pioneers of the theory was Imre Simon \cite {S}, a mathematician and computer
scientist from Brazil --- which is by the way the only reason for the peculiar
name ``tropical geometry''. Originally, the theory was developed in an applied
context of discrete mathematics and optimization, but it has not been part of
the mainstream in either of mathematics, computer science or engineering. Only
in the last few years have people realized its power for applications in fields
such as combinatorics, computational algebra, and algebraic geometry.

This is also why the theory of tropical algebraic geometry is still very much
in its beginnings: not even the concept of a variety has been defined yet in
tropical geometry in a general and satisfactory way. On the other hand there
are already many results in tropical geometry that show the power of these new
methods. For example, Mikhalkin has proven recently that tropical geometry can
be used to compute the numbers of plane curves of given genus $g$ and degree
$d$ through $ 3d+g-1 $ general points \cite {M} --- a deep result that had been
obtained first by Caporaso and Harris about ten years ago by a complicated
study of moduli spaces of plane curves \cite {CH}.

In this expository article we will for simplicity restrict ourselves mainly to
the well-established theory of tropical plane curves. Even in this special case
there are several seemingly different approaches to the theory. We will
describe these approaches in turn in chapter \ref {tropcurve} and discuss
possible generalizations at the end. We will then explain in chapter \ref
{classical} how some well-known results from classical geometry --- e.g.\ the
degree-genus formula and B\'ezout's theorem --- can be recovered (and reproven)
in the language of tropical geometry. Finally, in chapter \ref {enum} we will
discuss the most powerful applications of tropical geometry known so far,
namely to complex and real enumerative geometry.

%% file: tropcurve.tex
\section {Plane tropical curves} \label {tropcurve}


\subsection {Tropical curves as limits of amoebas} \label {amoebas}

With classical (complex) algebraic geometry in mind the most straightforward
way to tropical geometry is via so-called \emph {amoebas} of algebraic
varieties. For a complex plane curve $C$ the idea is simply to restrict it to
the open subset $ (\CC^*)^2 $ of the (affine or projective) plane and then to
map it to the \emph {real} plane by the map
\begin {align*}
  \Log: \qquad \quad (\CC^*)^2 &\to \RR^2 \\
    z=(z_1,z_2) &\mapsto (x_1,x_2):=(\log |z_1|, \log |z_2|).
\end {align*}
The resulting subset $ A = \Log (C \cap (\CC^*)^2) $ of $ \RR^2 $ is called the
\emph {amoeba} of the given curve. It is of course a two-dimensional subset of
$ \RR^2 $ since complex curves are real two-dimensional. The following picture
shows three examples (where we set $ e := \exp(1) $):

\begin {figure}[H]
  \begin {center} \input {pics/amoebas} \end {center}
  \caption {Three amoebas of plane curves}
  \label {pic-amoebas}
\end {figure}
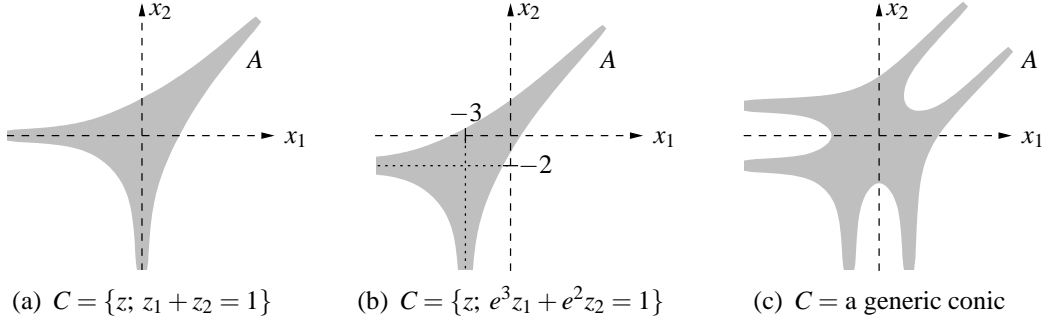

In fact, the shape of these pictures (that also explains the name ``amoeba'')
can easily be explained. In case (a) for example the curve $C$ contains exactly
one point whose $ z_1 $-coordinate is zero, namely $ (0,1) $. As $ \log 0 =
-\infty $ a small neighborhood of this point is mapped by $ \Log $ to the
``tentacle'' of the amoeba $A$ pointing to the left. In the same way a
neighborhood of $ (1,0) \in C $ leads to the tentacle pointing down, and points
of the form $ (z,1-z) $ with $ |z| \to \infty $ to the tentacle pointing to the
upper right.

In case (b) the multiplicative change in the variable simply leads to an
(additive) shift of the amoeba. In (c) a generic conic, i.e.\ a curve given by
a general polynomial of degree 2, has two points each where it meets the
coordinate axes, leading to two tentacles in each of the three directions. In
the same way one could consider curves of an arbitrary degree $d$ that would
give us amoebas with $d$ tentacles in each direction.

To make these amoebas into combinatorial objects the idea is simply to shrink
them to ``zero width''. So instead of the map $ \Log $ above let us consider
the maps
\begin {align*}
  \Log_t: \quad (\CC^*)^2 &\to \RR^2 \\
    (z_1,z_2) &\mapsto \left( - \log_t |z_1|, - \log_t |z_2| \right) =
      \left( -\frac {\log |z_1|}{\log t}, -\frac {\log |z_2|}{\log t} \right)
\end {align*}
for small $ t \in \RR $ and study the limit of the amoebas $ \Log_t (C \cap
(\CC^*)^2 ) $ as $t$ tends to zero. As $ \Log_t $ differs from $ \Log $ only by
a rescaling of the two axes the result for the curve in figure \ref
{pic-amoebas} (a) is the graph $ \Gamma $ shown in the following picture on the
left. We call $ \Gamma $ the \emph {tropical curve} determined by $C$.

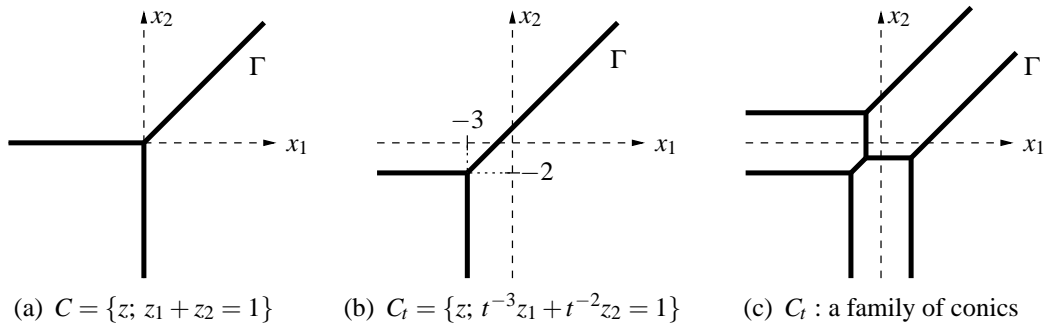
\begin {figure}[H]
  \begin {center} \input {pics/limits} \end {center}
  \caption {The tropical curves corresponding to the amoebas in figure \ref
    {pic-amoebas}}
  \label {pic-limits}
\end {figure}

If we did the same thing with the curve (b) the result would of course be that
we not only shrink the amoeba to zero width, but also move its ``vertex'' to
the origin, leading to the same tropical curve as in (a). To avoid this we
consider not only one curve $ C=\{z;\; e^3z_1+e^2z_2=1 \} $ but the \emph
{family} of curves $ C_t = \{ z;\; t^{-3}z_1+t^{-2}z_2=1 \} $ for small $ t \in
\RR $. This family has the property that $ C_t $ passes through $ (0,t^2) $ and
$ (t^3,0) $ for all $t$, and hence all $ \Log_t (C_t \cap (\CC^*)^2) $ have
their horizontal and vertical tentacles at $ z_2 = -2 $ and $ z_1 = -3 $,
respectively. So if we now take the limit as $ t \to 0 $ we shrink the width of
the amoeba to zero but keep its position in the plane: we get the shifted
tropical curve as in figure \ref {pic-limits} (b). We call this the tropical
curve determined by the family $ (C_t) $.

For (c) we can proceed in the same way: using a suitable family of conics we
can shrink the width of the amoeba $A$ to zero while keeping the position of
its tentacles fixed. The resulting tropical curve $ \Gamma $, i.e.\ the limit
of $ \Log_t (C_t \cap (\CC^*)^2) $, may e.g.\ look like in figure \ref
{pic-limits} (c). We will see at the end of section \ref {graphs} however
that this is not the only type of graph that we can obtain by a family of
conics in this way.

Summarizing we can say informally that a tropical curve should be a subset of $
\RR^2 $ obtained as the limit (in a certain sense) of the amoebas $ \Log_t (C_t
\cap (\CC^*)^2) $, where $ (C_t) $ is a suitable family of plane algebraic
curves. They are all piecewise linear graphs with certain properties that we
will study later.

\subsection {Tropical curves via varieties over the field of Puiseux series}
  \label {puiseux}

Of course this method of constructing (or even defining) tropical curves is
very cumbersome as it always involves a limiting process over a whole family of
complex curves. There is an elegant way to hide this limiting process by
replacing the ground field $ \CC $ by the field $K$ of so-called \emph {Puiseux
series}, i.e.\ the field of formal power series $ a = \sum_{q \in \QQ} a_q t^q
$ in a variable $t$ such that the subset of $ \QQ $ of all $q$ with $ a_q \neq
0 $ is bounded below and has a finite set of denominators. For such an $ a \in
K $ with $ a \neq 0 $ the infimum of all $q$ with $ a_q \neq 0 $ is actually a
minimum; it is called the \emph {valuation} of $a$ and denoted $ \val a $.

Using this construction we can say for example that our family (b) in section
\ref {amoebas}
  \[ C = \{ z \in K^2;\; t^{-3}z_1+t^{-2}z_2 = 1 \} \subset K^2 \]
now defines \emph {one single curve} in the affine plane over this new field
$K$. How do we now perform the limit $ t \to 0 $ in this set-up? For an
element $ a = \sum_{q \in \QQ} a_q t^q \in K $ only the term with the smallest
exponent, i.e.\ $ a_{\val a} t^{\val a} $, will be relevant in this limit. So
applying the map $ \log_t $ we get for small $t$
  \[ \log_t |a| \approx \log_t |a_{\val a} t^{\val a}| = \val a +
       \log_t |a_{\val a}| \approx \val a. \]
In our new picture the operations of applying the map $ \Log_t $ and taking the
limit for $ t \to 0 $ therefore correspond to the map
\begin {align*}
  \Val: (K^*)^2 &\to \RR^2 \\
    (z_1,z_2) &\mapsto (x_1,x_2) := (-\val z_1, -\val z_2).
\end {align*}
Using this observation we can now give our first rigorous definition of plane
tropical curves:

\textbf {Definition A.} \emph {A plane tropical curve is a subset of $ \RR^2 $
of the form $ \Val (C \cap (K^*)^2) $, where $C$ is a plane algebraic curve in
$ K^2 $. (Strictly speaking we should take the closure of $ \Val (C \cap
(K^*)^2) $ in $ \RR^2 $ since the image of the valuation map $ \Val $ is by
definition contained in $ \QQ^2 $.)}

Note that this definition is now purely algebraic and does not involve any
limit taking processes. As $K$ is an algebraically closed field of
characteristic zero (in fact it is the algebraic closure of the field of
Laurent series in $t$) the theory of algebraic geometry of plane curves over
$K$ is largely identical to that of algebraic curves over $ \CC $.

As an example let us consider again case (b) of section \ref {amoebas}, i.e.\
the curve $ C \subset K^2 $ given by the equation $ t^{-3}z_1+t^{-2}z_2 = 1 $.
If $ (z_1,z_2) \in C \cap (K^*)^2 $ then $ \Val (z_1,z_2) $ can give three
different kinds of results:
\begin {itemize}
\item If $ \val z_1 > 3 $ then the valuation of $ z_2 = t^2 - t^{-1} z_1 $ is 2
  since all exponents of $t$ in $ t^{-1} z_1 $ are bigger than 2. Hence these
  points map precisely to the left edge of the tropical curve in figure \ref
  {pic-limits} (b) under $ \Val $.
\item In the same way we get the bottom edge of this tropical curve if $ \val
  z_2 > 2 $.
\item If $ \val z_1 \le 3 $ and $ \val z_2 \le 2 $ then the equation $
  t^{-3}z_1+t^{-2}z_2 = 1 $ shows that the leading terms of $ t^{-3}z_1 $ and $
  t^{-2}z_2 $ must have the same valuation, i.e.\ that $ \val z_1 = \val z_2 +
  1 $. This leads to the upper right edge of the tropical curve in figure \ref
  {pic-limits} (b).
\end {itemize}
So we recover our old result, i.e.\ the tropical curve drawn in figure \ref
{pic-limits} (b).

One special case is worth mentioning: if the curve $ C \subset K^2 $ is given
by an equation whose coefficients lie in $ \CC $ (i.e.\ are ``independent of
$t$'') then for any point $ (z_1(t),z_2(t)) \in C $ the points $
(z_1(t^q),z_2(t^q)) $ for $ q \in \QQ $ are obviously in $C$ as well. As
replacing $t$ by $ t^q $ for some $ q>0 $ simply multiplies the valuation with
$q$ we conclude that the tropical variety associated to $C$ in this case is a
\emph {cone} (i.e.\ a union of half-rays starting at the origin) --- as it was
the case e.g.\ in figure \ref {pic-limits} (a).

\subsection {Tropical curves as varieties over the max-plus semiring}
  \label {maxplus}

We now want to study definition A in more detail. Let $ C \subset K^2 $ be a
plane algebraic curve given by the polynomial equation
  \[ C = \left\{ (z_1,z_2) \in K^2 ;\;
         f(z_1,z_2) := \sum_{i,j \in \NN} a_{ij} z_1^i z_2^j = 0 \right\} \]
for some $ a_{ij} \in K $ of which only finitely many are non-zero. Note that
the valuation of a summand of $ f(z_1,z_2) $ is
  \[ \val (a_{ij} z_1^i z_2^j) = \val a_{ij} + i \val z_1 + j \val z_2. \]
Now if $ (z_1,z_2) $ is a point of $C$ then all these summands add up to zero.
In particular, the lowest valuation of these summands must occur at least twice
since otherwise the corresponding terms in the sum could not cancel. For the
corresponding point $ (x_1,x_2) = \Val (z_1,z_2) = (-\val z_1,-\val z_2) $ of
the tropical curve this obviously means that in the expression
  \[ g(x_1,x_2) := \max \{ ix_1 + j x_2 - \val a_{ij} ;\;
       (i,j) \in \NN^2 \mbox { with $ a_{ij} \neq 0 $} \} \tag {$*$} \]
the maximum is taken on at least twice. It follows that the tropical curve
determined by $C$ is contained in the ``corner locus'' of this convex piecewise
linear function $g$, i.e.\ in the locus where this function is not
differentiable. In fact, Kapranov's theorem states that the converse inclusion
holds as well, i.e.\ that the tropical curve determined by $C$ is precisely
this corner locus (see e.g.\ \cite {K}, \cite {Sh}).

As an example let us consider again the curve $ C=\{ (z_1,z_2) ;\; t^{-3} z_1 +
t^{-2} z_2 - 1 = 0 \} \subset K^2 $ that we discussed in the previous section.
The corresponding convex piecewise linear function is
  \[ g(x_1,x_2) = \max \{ x_1+3, x_2+2,0 \} \]
The following picture shows how the corner locus of this function gives us back
the tropical curve of figure \ref {pic-limits} (b).

\begin {figure}[H]
  \begin {center} \input {pics/corner} \end {center}
  \caption {A tropical curve as the corner locus of a convex piecewise linear
    function}
  \label {pic-corner}
\end {figure}
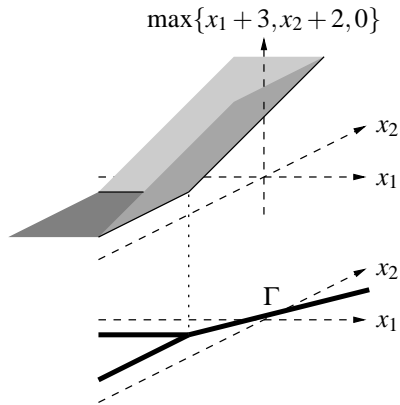

These convex piecewise linear functions are often written in a different way in
order to resemble the notation of the original polynomial: for two real numbers
$ x,y $ we define ``tropical addition'' and ``tropical multiplication'' simply
by
  \[ x \oplus y := \max \{x,y\}
       \qquad \mbox {and} \qquad
     x \odot y := x + y. \]
The real numbers together with these two operations form a semiring, i.e.\ they
satisfy all properties of a ring except for the existence of additive neutral
and inverse elements. Sometimes an element $ -\infty $ is formally added to the
real numbers to serve as a neutral element, but there is certainly no way to
construct inverse elements as this would require equations of the form $ \max
\{ -\infty,x\} = -\infty $ to be solvable.

Using this notation we can write our convex piecewise linear function $(*)$
above as
  \[ g(x_1,x_2) = \bigoplus_{i,j} (-\val a_{ij}) \odot x_1^{\odot i} \odot
    x_2^{\odot j}. \]
We call this expression the \emph {tropicalization} of the original polynomial
$f$. It can be considered as a ``tropical polynomial'', i.e.\ as a polynomial
in the tropical semiring. For example, the tropicalization of the polynomial $
t^{-3} z_1 + t^{-2} z_2 - 1 $ is just
  \[ 3 \odot x_1 \oplus 2 \odot x_2 \oplus 0 = \max \{x_1+3,x_2+2,0\}. \]
(Note that the addition of 0 is not superfluous here since 0 is not a neutral
element for tropical addition!)

We can therefore now give an alternative definition of plane tropical curves
that does not involve the somewhat complicated field of Puiseux series any
more:

\textbf {Definition B.} \emph {A plane tropical curve is a subset of $ \RR^2 $
that is the corner locus of a tropical polynomial, i.e.\ of a polynomial in the
tropical semiring $ (\RR,\oplus,\odot) = (\RR,\max,+) $.}

Again there is a special case that is completely analogous to the one mentioned
at the end of section \ref {puiseux}: if the tropical polynomial $g$ is the
maximum of linear functions without constant terms (e.g.\ because it is the
tropicalization of a polynomial with coefficients that do not depend on $t$)
then the corner locus of $g$ is a cone. If this is not the case and $g$ is the
maximum of many affine functions then its corner locus will in general be a
complicated piecewise linear graph in the plane as e.g.\ in figure \ref
{pic-limits} (c).

\subsection {Tropical curves as balanced graphs} \label {graphs}

Our definition B now allows us to give an easy and entirely geometric
characterization of plane tropical curves. We have already seen that a tropical
curve $ \Gamma $ is a graph in $ \RR^2 $ whose edges are line segments. Let us
consider $ \Gamma $ locally around a vertex $ V \in \Gamma $. For simplicity we
shift coordinates so that $V$ is the origin in $ \RR^2 $ and thus $ \Gamma $
becomes a cone locally around $V$. We have seen already at the end of section
\ref {maxplus} that $ \Gamma $ is then locally the corner locus of a tropical
polynomial of the form
  \[ g(x_1,x_2) = \bigoplus_i \;
     x_1^{\odot a_1^{(i)}} \odot x_2^{\odot a_2^{(i)}}
       = \max \{ a_1^{(i)} x_1 + a_2^{(i)} x_2 ;\; i=1,\dots,n \}. \]
for some $ a^{(i)} = (a_1^{(i)},a_2^{(i)}) \in \NN^2 $. Let $ \Delta $ be the
convex hull of the points $ a^{(i)} $, as indicated in the following picture
on the left:

\begin {figure}[H]
  \begin {center} \input {pics/balance} \end {center}
  \caption {A local picture of a plane tropical curve}
  \label {pic-balance}
\end {figure}
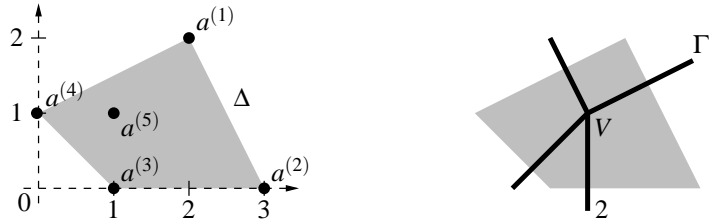

First of all we claim that any point $ a^{(i)} $ that is not a vertex of $
\Delta $ is irrelevant for the tropical curve $ \Gamma $. In fact, it is
impossible for such an $ a^{(i)} $ (as e.g.\ $ a^{(5)} $ in the example above)
that the expression $ a_1^{(i)} x_1 + a_2^{(i)} x_2 $ is strictly bigger than
all the other $ a_1^{(j)} x_1 + a_2^{(j)} x_2 $ for some $ x_1,x_2 \in \RR $.
Hence $g$ and therefore also its corner locus remain the same if we drop this
term. In particular we see that --- unlike in classical algebraic geometry ---
\emph {there is no hope for a one-to-one correspondence between tropical curves
and tropical polynomials (up to scalars)}.

It is now easy to see that the corner locus of $g$ consists precisely of those
points where
  \[ g(x_1,x_2) = a_1^{(i)} x_1 + a_2^{(i)} x_2 = a_1^{(j)} x_1 + a_2^{(j)} x_2
     \]
for two adjacent vertices $ a^{(i)} $ and $ a^{(j)} $ of $ \Delta $. An easy
computation shows that for fixed $i$ and $j$ this is precisely the half-ray
starting from the origin and pointing in the direction of the outward normal
of the edge joining $ a^{(i)} $ and $ a^{(j)} $. So as shown in figure \ref
{pic-balance} on the right the tropical curve $ \Gamma $ is simply the union of
all these outward normal lines locally around $V$. In particular all edges of $
\Gamma $ have rational slopes.

There is one more important condition on the edges of $ \Gamma $ around $V$
that follows from this observation. If $ a^{(1)},\dots,a^{(n)} $ are the
vertices of $ \Delta $ in clockwise direction then an outward normal vector of
the edge joining $ a^{(i)} $ and $ a^{(i+1)} $ (where we set $ a^{(n+1)} :=
a^{(1)} $) is $ v^{(i)} := (a_2^{(i)}-a_2^{(i+1)},a_1^{(i+1)}-a_1^{(i)}) $ for
all $ i=1,\dots,n $. In particular it follows that $ \sum_{i=1}^n v^{(i)} = 0
$. This fact is usually expressed as follows: we write the vectors $ v^{(i)} $
as $ v^{(i)} = w^{(i)} \cdot u^{(i)} $ where $ u^{(i)} $ is the primitive
integral vector in the direction of $ v^{(i)} $ and $ w^{(i)} \in \NN_{>0} $.
We call $ w^{(i)} $ the \emph {weight} of the corresponding edge of $ \Gamma $
and thus consider $ \Gamma $ to be a weighted graph. Our equation $ \sum_i
v^{(i)} = 0 $ then states that the weighted sum of the primitive integral
vectors of the edges around every vertex of $ \Gamma $ is 0. This is usually
called the \emph {balancing condition}. For example, in figure \ref
{pic-balance} the edge of $ \Gamma $ pointing down has weight 2 (since $
v^{(2)} = (0,-2) = 2 \cdot (0,-1) $), whereas all other edges have weight 1.
In this paper we will usually label the edges with their corresponding weights
unless these weights are 1. The balancing condition around the vertex $V$ then
reads
  \[ (2,1) + 2 \cdot (0,-1) + (-1,-1) + (-1,2) = (0,0) \]
in this example.

Together with our observation of section \ref {amoebas} that (at least generic)
plane algebraic curves of degree $d$ lead to plane tropical curves with $d$
ends each in the directions $ (-1,0) $, $ (0,-1) $, and $ (1,1) $, we arrive at
the following somewhat longer but purely geometric definition of plane tropical
curves:

\textbf {Definition C.} \emph {A plane tropical curve of degree $d$ is a
weighted graph $ \Gamma $ in $ \RR^2 $ such that
\begin {enumerate}
\item every edge of $ \Gamma $ is a line segment with rational slope;
\item $ \Gamma $ has $d$ ends each in the directions $ (-1,0) $, $ (0,-1) $,
  and $ (1,1) $ (where an end of weight $w$ counts $w$ times);
\item at every vertex $V$ of $ \Gamma $ the balancing condition holds: the
  weighted sum of the primitive integral vectors of the edges around $V$ is
  zero.
\end {enumerate}
}%

Strictly speaking we have only explained above why a plane tropical curve in
the sense of definition B gives rise to a curve in the sense of definition C.
One can show that the converse holds as well; a proof can e.g.\ be found in
\cite {M} or \cite {Sp} chapter 5.

With this definition it has now become a combinatorial problem to find all
types of plane tropical curves of a given degree. In fact, the construction
given above globalizes well: assume that $ \Gamma $ is the tropical curve given
as the corner locus of the tropical polynomial
  \[ g(x_1,x_2) = \max \{ a_1^{(i)} x_1 + a_2^{(i)} x_2 + b^{(i)} ;\;
       i=1,\dots,n \}. \]
If $g$ is the tropicalization of a polynomial of degree $d$ then the $ a^{(i)}
$ are all integer points in the triangle $ \Delta_d := \{ (a_1,a_2);\; a_1 \ge
0, a_2 \ge 0, a_1+a_2 \le d \} $. Consider two terms $ i,j \in \{1,\dots,n\} $
with $ a^{(i)} \neq a^{(j)} $. If there is a point $ (x_1,x_2) \in \RR^2 $
such that
  \[ g(x_1,x_2) = a_1^{(i)} x_1 + a_2^{(i)} x_2 + b^{(i)}
                = a_1^{(j)} x_1 + a_2^{(j)} x_2 + b^{(j)}
     \]
then we draw a straight line in $ \Delta_d $ through the points $ a^{(i)} $ and
$ a^{(j)} $. This way we obtain a subdivision of $ \Delta_d $ whose edges
correspond to the edges of $ \Gamma $ and whose 2-dimensional cells correspond
to the vertices of $ \Gamma $. This subdivision is usually called the \emph
{Newton subdivision} corresponding to $ \Gamma $. So to find all types of plane
tropical curves of degree $d$ one has to list all subdivisions of $ \Delta_d $
and check which of them are induced by a tropical curve as above. As the
simplest example there is only the trivial subdivision of $ \Delta_1 $, leading
to the only type of plane tropical curve of degree 1:

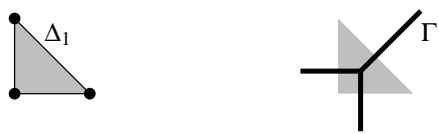
\begin {figure}[H]
  \begin {center} \input {pics/lines} \end {center}
  \caption {Tropical lines are ``dual'' to $ \Delta_1 $}
  \label {pic-lines}
\end {figure}

The following picture shows all non-degenerated cases in degree 2, where by
``non-degenerate'' we mean that the subdivision of $ \Delta_2 $ is maximal.

\begin {figure}[H]
  \begin {center} \input {pics/conics} \end {center}
  \caption {The four types of (non-degenerated) tropical plane conics}
  \label {pic-conics}
\end {figure}
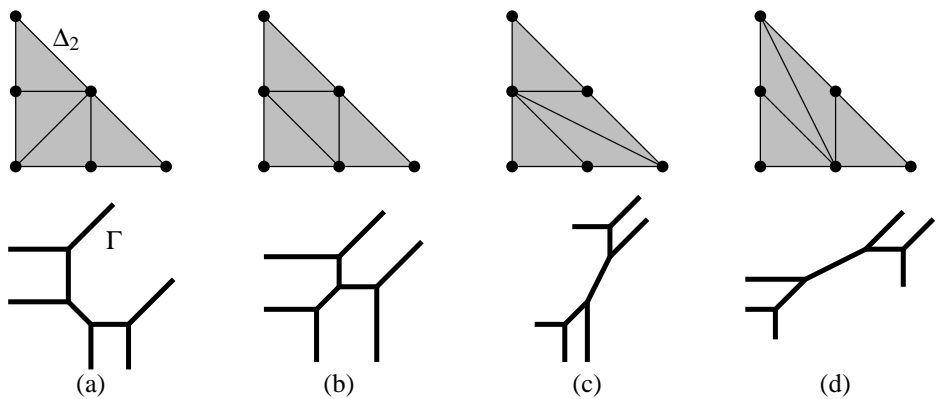

In general, the problem of finding an algorithm to generate all such
subdivisions has already been studied extensively in geometric combinatorics
\cite {IMTI,R}.

Note however that not every subdivision gives rise to a type of tropical
curves. First of all it is obvious by our constructions above that we need
subdivisions into \emph {convex} polytopes in order to have a tropical curve
corresponding to it. But this condition is not sufficient: assume for example
that we have the following (local) picture $ \Delta $ somewhere in the Newton
subdivision of a plane tropical curve:

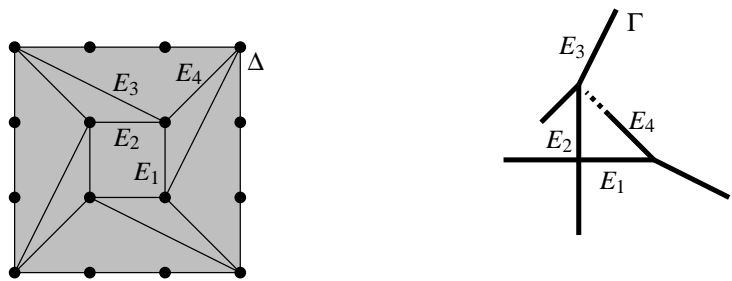
\begin {figure}[H]
  \begin {center} \input {pics/convex} \end {center}
  \caption {A subdivision that is not induced by a tropical curve}
  \label {pic-convex}
\end {figure}

In the tropical curve $ \Gamma $ that would correspond to this subdivision
the edge $ E_4 $ would have to meet $ E_3 $ and not $ E_2 $ at the dotted end,
so $ E_2 $ must be longer than $ E_1 $. But the same argument can be used
cyclically to conclude that each edge around the central vertex must be longer
than the previous one. As this is not possible we conclude that there cannot be
a tropical curve corresponding to this subdivision of $ \Delta $. In fact, the
subdivisions corresponding to tropical curves are precisely the ones that are
usually called the \emph {regular polyhedral subdivisions}, i.e.\ the ones that
can be written as the corner locus of a piecewise linear convex function
similarly to figure \ref {pic-corner}.

It should also be stressed that the subdivision of $ \Delta_d $ determines only
the combinatorial type of the tropical curve and not the curve itself. For
example, each of the types in figure \ref {pic-conics} describes a (real)
5-dimensional family of plane tropical conics since the lengths of the bounded
edges (3 parameters) as well as the position in the plane (2 parameters) can
vary arbitrarily. Note that this agrees nicely with the classical picture:
conics in the complex plane vary in a 5-dimensional family as well
(corresponding to the 6 coefficients of a quadratic equation modulo a common
scalar).

\subsection {Generalizations} \label {generalizations}

At the end of this chapter let us brief\/ly describe how the theory of plane
tropical curves given above can be generalized.

First of all it is quite obvious to note that essentially the same
constructions and the same theory can be carried through for curves that are
not necessarily in the plane but in any \emph {toric surface}, i.e.\ in any
surface $X$ with a $ \CC^* $-action that contains $ (\CC)^* $ as a dense open
subset. Definition A remains unchanged in this case; the resulting tropical
curves will still be graphs in $ \RR^2 $. Definition B only has to be modified
as to allow Laurent polynomials compatible with the chosen homology class of
the curves; and in definition C the only change is in the directions of the
ends of $ \Gamma $. In fact, a toric surface $X$ together with a positive
homology class corresponds exactly to a convex integral polytope $ \Delta $,
and the types of tropical curves coming from this homology class are precisely
those determined by subdivisions of $ \Delta $ as explained at the end of
section \ref {graphs}.

In the same way one can also consider general hypersurfaces instead of plane
curves. Except for the existence of more variables there are no changes in
definitions A and B, and there is an analogous version of the balancing
condition and definition C too. The resulting tropical hypersurfaces are
weighted polyhedral complexes in a real vector space. It is also true in this
case that the combinatorial types of hypersurfaces correspond to subdivisions
of a higher-dimensional polytope.

The theory becomes more difficult however in the case of varieties of higher
codimension, e.g.\ space curves. It is probably agreed upon that definition A
would be the ``correct'' one also in this case, i.e.\ that tropical varieties
are by definition the images of classical varieties over the field of Puiseux
series under the valuation map. This definition is however hard to work with in
practice --- it would be much more convenient to think of tropical
hypersurfaces as in definition B and of general tropical varieties as
intersections of such tropical hypersurfaces. Unfortunately, if $ X \subset K^n
$ is a variety given by some polynomial equations $ f_1 = \cdots = f_r = 0 $ it
is (in contrast to the codimension-1 case) in general \emph {not} true that the
tropical variety corresponding to $X$ (i.e.\ the image of $X$ under the
valuation map) is simply given by the intersection of the corner loci of the
tropicalizations of $ f_1,\dots,f_r $. In fact, it is not even true in general
that the intersection of tropical hypersurfaces is a tropical variety at all:
if we intersect e.g.\ a tropical line as in figure \ref {pic-limits} (b) with
the same line shifted a bit to the left then the result is a single half-ray
--- which is not a tropical variety. It has been shown however that the
tropical variety corresponding to $X$ can always be written as an intersection
of the corner loci of the tropicalizations of (finitely many) suitably chosen
generators of the ideal $ (f_1,\dots,f_r) $ defining $X$ \cite {BJSST}.
Moreover, there exists an implemented algorithm to perform these computations
explicitly \cite {BJSST,J} so that --- from an algorithmic point of view ---
"`every variety can be tropicalized"'. However, the necessary calculations rely
on Gr\"obner basis techniques and thus become complicated very soon (more
precisely, most questions regarding the computation of tropical varieties are
NP-hard in the language of complexity analysis \cite {T}).

As for definition C the conditions listed there can be adapted to make sense in
higher codimensions as well, so one could try to use definition C and say that
e.g.\ space curves are simply balanced graphs in $ \RR^3 $. Unfortunately it
turns out that this definition would not be equivalent to definition A. While
it is true that every tropical space curve in the sense of definition A gives
rise to a balanced graph in $ \RR^3 $ in the sense of definition C the converse
does not hold in general \cite {Sp}. There is no known general criterion yet
to decide exactly when a balanced graph in $ \RR^3 $ can be obtained as the
image of a curve in $ K^3 $ under the valuation map, although a sufficient
criterion is given in \cite {Sp}.

Finally one should note that even the most generally applicable definition A
depends on a given embedding of the original variety over $K$ in an affine or
projective space (or a toric variety). There is no (known) way to associate a
tropical variety to a given abstract variety over $K$. In fact, there is not
even a good theory yet of what an abstract tropical variety should be. There is
some recent work of Mikhalkin however that tries to build up a theory of
tropical geometry completely in parallel to algebraic geometry, replacing the
ground field by the semiring $ (\RR,\oplus,\odot) $ \cite {M2}.

%% file: pics/amoebas.tex
\begin{picture}(0,0)%
\includegraphics{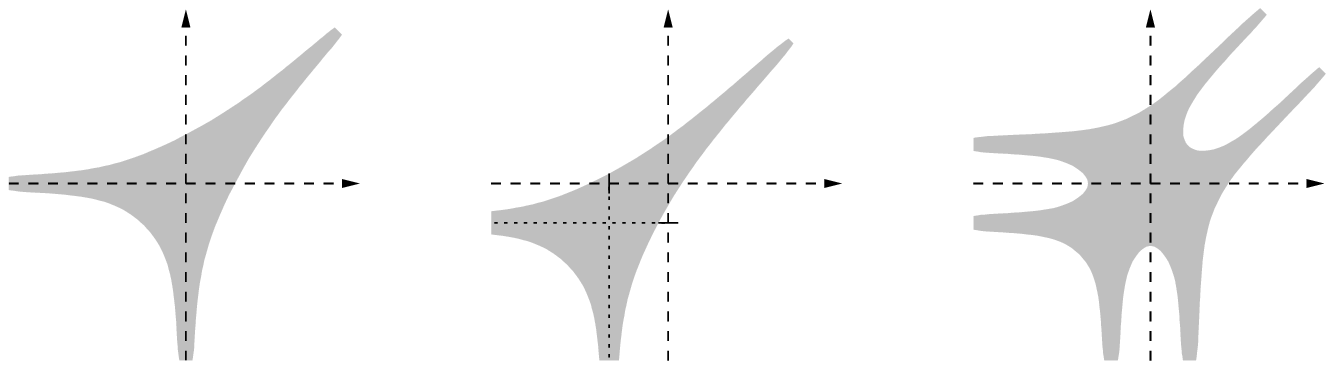}%
\end{picture}%
\setlength{\unitlength}{4144sp}%
\begingroup\makeatletter\ifx\SetFigFont\undefined%
\gdef\SetFigFont#1#2#3#4#5{%
  \reset@font\fontsize{#1}{#2pt}%
  \fontfamily{#3}\fontseries{#4}\fontshape{#5}%
  \selectfont}%
\fi\endgroup%
\begin{picture}(6257,2045)(675,-1610)
\put(2431,-556){\makebox(0,0)[lb]{\smash{{\SetFigFont{10}{12.0}{\familydefault}{\mddefault}{\updefault}{\color[rgb]{0,0,0}$x_1$}%
}}}}
\put(1621,209){\makebox(0,0)[lb]{\smash{{\SetFigFont{10}{12.0}{\familydefault}{\mddefault}{\updefault}{\color[rgb]{0,0,0}$x_2$}%
}}}}
\put(4636,-556){\makebox(0,0)[lb]{\smash{{\SetFigFont{10}{12.0}{\familydefault}{\mddefault}{\updefault}{\color[rgb]{0,0,0}$x_1$}%
}}}}
\put(6841,-556){\makebox(0,0)[lb]{\smash{{\SetFigFont{10}{12.0}{\familydefault}{\mddefault}{\updefault}{\color[rgb]{0,0,0}$x_1$}%
}}}}
\put(6031,209){\makebox(0,0)[lb]{\smash{{\SetFigFont{10}{12.0}{\familydefault}{\mddefault}{\updefault}{\color[rgb]{0,0,0}$x_2$}%
}}}}
\put(3826,209){\makebox(0,0)[lb]{\smash{{\SetFigFont{10}{12.0}{\familydefault}{\mddefault}{\updefault}{\color[rgb]{0,0,0}$x_2$}%
}}}}
\put(5986,-1546){\makebox(0,0)[b]{\smash{{\SetFigFont{10}{12.0}{\familydefault}{\mddefault}{\updefault}{\color[rgb]{0,0,0}(c) $ \; C = $ a generic conic}%
}}}}
\put(2206,-106){\makebox(0,0)[lb]{\smash{{\SetFigFont{10}{12.0}{\familydefault}{\mddefault}{\updefault}{\color[rgb]{0,0,0}$A$}%
}}}}
\put(6841,-106){\makebox(0,0)[lb]{\smash{{\SetFigFont{10}{12.0}{\familydefault}{\mddefault}{\updefault}{\color[rgb]{0,0,0}$A$}%
}}}}
\put(1576,-1546){\makebox(0,0)[b]{\smash{{\SetFigFont{10}{12.0}{\familydefault}{\mddefault}{\updefault}{\color[rgb]{0,0,0}(a) $ \; C=\{ z;\;z_1+z_2=1 \} $}%
}}}}
\put(3826,-736){\makebox(0,0)[lb]{\smash{{\SetFigFont{10}{12.0}{\familydefault}{\mddefault}{\updefault}{\color[rgb]{0,0,0}$-2$}%
}}}}
\put(3511,-421){\makebox(0,0)[b]{\smash{{\SetFigFont{10}{12.0}{\familydefault}{\mddefault}{\updefault}{\color[rgb]{0,0,0}$-3$}%
}}}}
\put(4321,-106){\makebox(0,0)[lb]{\smash{{\SetFigFont{10}{12.0}{\familydefault}{\mddefault}{\updefault}{\color[rgb]{0,0,0}$A$}%
}}}}
\put(3781,-1546){\makebox(0,0)[b]{\smash{{\SetFigFont{10}{12.0}{\familydefault}{\mddefault}{\updefault}{\color[rgb]{0,0,0}(b) $ \; C=\{z;\;e^3z_1+e^2z_2=1\} $}%
}}}}
\end{picture}%

%% file: pics/limits.tex
\begin{picture}(0,0)%
\includegraphics{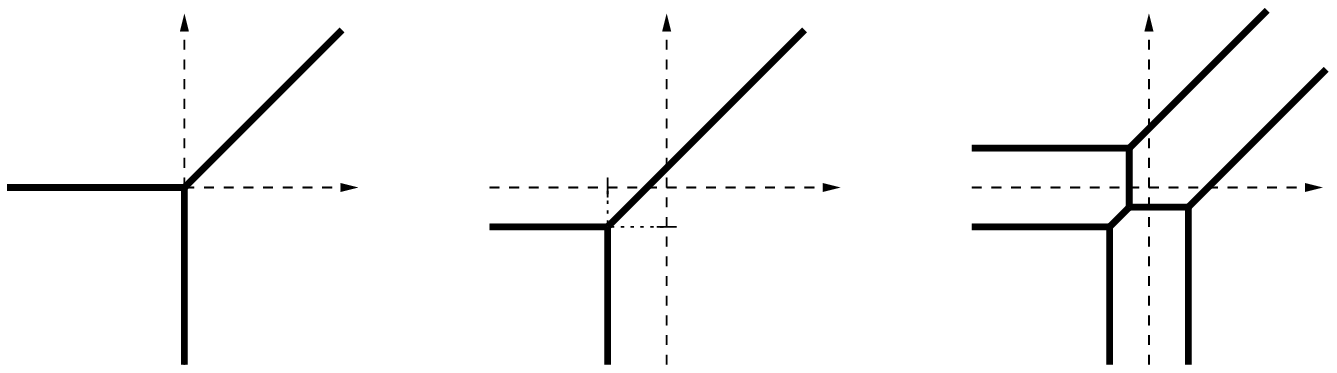}%
\end{picture}%
\setlength{\unitlength}{4144sp}%
\begingroup\makeatletter\ifx\SetFigFont\undefined%
\gdef\SetFigFont#1#2#3#4#5{%
  \reset@font\fontsize{#1}{#2pt}%
  \fontfamily{#3}\fontseries{#4}\fontshape{#5}%
  \selectfont}%
\fi\endgroup%
\begin{picture}(6123,1966)(733,-1610)
\put(2431,-556){\makebox(0,0)[lb]{\smash{{\SetFigFont{10}{12.0}{\familydefault}{\mddefault}{\updefault}{\color[rgb]{0,0,0}$x_1$}%
}}}}
\put(1621,209){\makebox(0,0)[lb]{\smash{{\SetFigFont{10}{12.0}{\familydefault}{\mddefault}{\updefault}{\color[rgb]{0,0,0}$x_2$}%
}}}}
\put(4636,-556){\makebox(0,0)[lb]{\smash{{\SetFigFont{10}{12.0}{\familydefault}{\mddefault}{\updefault}{\color[rgb]{0,0,0}$x_1$}%
}}}}
\put(6841,-556){\makebox(0,0)[lb]{\smash{{\SetFigFont{10}{12.0}{\familydefault}{\mddefault}{\updefault}{\color[rgb]{0,0,0}$x_1$}%
}}}}
\put(6031,209){\makebox(0,0)[lb]{\smash{{\SetFigFont{10}{12.0}{\familydefault}{\mddefault}{\updefault}{\color[rgb]{0,0,0}$x_2$}%
}}}}
\put(3826,209){\makebox(0,0)[lb]{\smash{{\SetFigFont{10}{12.0}{\familydefault}{\mddefault}{\updefault}{\color[rgb]{0,0,0}$x_2$}%
}}}}
\put(1576,-1546){\makebox(0,0)[b]{\smash{{\SetFigFont{10}{12.0}{\familydefault}{\mddefault}{\updefault}{\color[rgb]{0,0,0}(a) $ \; C=\{ z;\;z_1+z_2=1 \} $}%
}}}}
\put(2206,-106){\makebox(0,0)[lb]{\smash{{\SetFigFont{10}{12.0}{\familydefault}{\mddefault}{\updefault}{\color[rgb]{0,0,0}$\Gamma$}%
}}}}
\put(6841,-106){\makebox(0,0)[lb]{\smash{{\SetFigFont{10}{12.0}{\familydefault}{\mddefault}{\updefault}{\color[rgb]{0,0,0}$\Gamma$}%
}}}}
\put(4366,-61){\makebox(0,0)[lb]{\smash{{\SetFigFont{10}{12.0}{\familydefault}{\mddefault}{\updefault}{\color[rgb]{0,0,0}$\Gamma$}%
}}}}
\put(3826,-736){\makebox(0,0)[lb]{\smash{{\SetFigFont{10}{12.0}{\familydefault}{\mddefault}{\updefault}{\color[rgb]{0,0,0}$-2$}%
}}}}
\put(3511,-421){\makebox(0,0)[b]{\smash{{\SetFigFont{10}{12.0}{\familydefault}{\mddefault}{\updefault}{\color[rgb]{0,0,0}$-3$}%
}}}}
\put(5986,-1546){\makebox(0,0)[b]{\smash{{\SetFigFont{10}{12.0}{\familydefault}{\mddefault}{\updefault}{\color[rgb]{0,0,0}(c) $ \; C_t: $ a family of conics}%
}}}}
\put(3781,-1546){\makebox(0,0)[b]{\smash{{\SetFigFont{10}{12.0}{\familydefault}{\mddefault}{\updefault}{\color[rgb]{0,0,0}(b) $ \; C_t=\{z;\;t^{-3}z_1+t^{-2}z_2=1\} $}%
}}}}
\end{picture}%

%% file: pics/corner.tex
\begin{picture}(0,0)%
\includegraphics{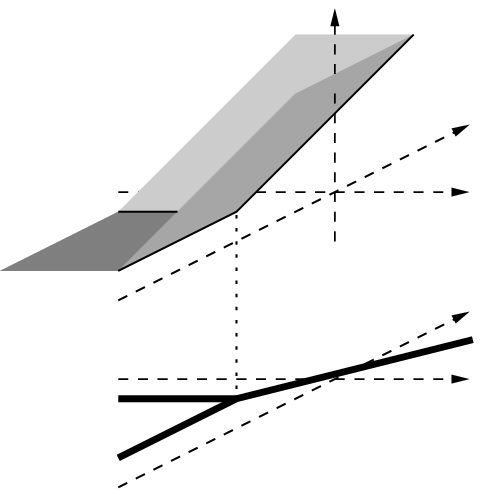}%
\end{picture}%
\setlength{\unitlength}{4144sp}%
\begingroup\makeatletter\ifx\SetFigFont\undefined%
\gdef\SetFigFont#1#2#3#4#5{%
  \reset@font\fontsize{#1}{#2pt}%
  \fontfamily{#3}\fontseries{#4}\fontshape{#5}%
  \selectfont}%
\fi\endgroup%
\begin{picture}(2221,2409)(1125,-2863)
\put(3331,-1546){\makebox(0,0)[lb]{\smash{{\SetFigFont{10}{12.0}{\familydefault}{\mddefault}{\updefault}{\color[rgb]{0,0,0}$ x_1 $}%
}}}}
\put(3331,-2401){\makebox(0,0)[lb]{\smash{{\SetFigFont{10}{12.0}{\familydefault}{\mddefault}{\updefault}{\color[rgb]{0,0,0}$ x_1 $}%
}}}}
\put(3331,-2086){\makebox(0,0)[lb]{\smash{{\SetFigFont{10}{12.0}{\familydefault}{\mddefault}{\updefault}{\color[rgb]{0,0,0}$ x_2 $}%
}}}}
\put(3331,-1231){\makebox(0,0)[lb]{\smash{{\SetFigFont{10}{12.0}{\familydefault}{\mddefault}{\updefault}{\color[rgb]{0,0,0}$ x_2 $}%
}}}}
\put(2656,-601){\makebox(0,0)[b]{\smash{{\SetFigFont{10}{12.0}{\familydefault}{\mddefault}{\updefault}{\color[rgb]{0,0,0}$ \max \{ x_1+3,x_2+2,0 \} $}%
}}}}
\put(2746,-2266){\makebox(0,0)[rb]{\smash{{\SetFigFont{10}{12.0}{\familydefault}{\mddefault}{\updefault}{\color[rgb]{0,0,0}$ \Gamma $}%
}}}}
\end{picture}%

%% file: pics/balance.tex
\begin{picture}(0,0)%
\includegraphics{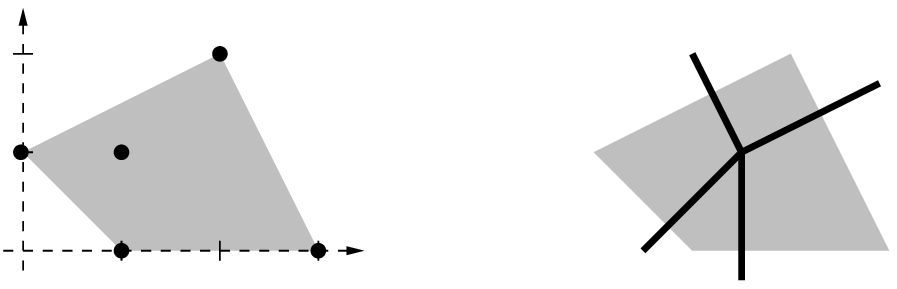}%
\end{picture}%
\setlength{\unitlength}{4144sp}%
\begingroup\makeatletter\ifx\SetFigFont\undefined%
\gdef\SetFigFont#1#2#3#4#5{%
  \reset@font\fontsize{#1}{#2pt}%
  \fontfamily{#3}\fontseries{#4}\fontshape{#5}%
  \selectfont}%
\fi\endgroup%
\begin{picture}(4066,1372)(346,-1646)
\put(901,-1591){\makebox(0,0)[b]{\smash{{\SetFigFont{10}{12.0}{\familydefault}{\mddefault}{\updefault}{\color[rgb]{0,0,0}$1$}%
}}}}
\put(1801,-1591){\makebox(0,0)[b]{\smash{{\SetFigFont{10}{12.0}{\familydefault}{\mddefault}{\updefault}{\color[rgb]{0,0,0}$3$}%
}}}}
\put(1351,-1591){\makebox(0,0)[b]{\smash{{\SetFigFont{10}{12.0}{\familydefault}{\mddefault}{\updefault}{\color[rgb]{0,0,0}$2$}%
}}}}
\put(361,-1006){\makebox(0,0)[rb]{\smash{{\SetFigFont{10}{12.0}{\familydefault}{\mddefault}{\updefault}{\color[rgb]{0,0,0}$1$}%
}}}}
\put(361,-556){\makebox(0,0)[rb]{\smash{{\SetFigFont{10}{12.0}{\familydefault}{\mddefault}{\updefault}{\color[rgb]{0,0,0}$2$}%
}}}}
\put(1621,-916){\makebox(0,0)[lb]{\smash{{\SetFigFont{10}{12.0}{\familydefault}{\mddefault}{\updefault}{\color[rgb]{0,0,0}$ \Delta $}%
}}}}
\put(1396,-466){\makebox(0,0)[lb]{\smash{{\SetFigFont{10}{12.0}{\familydefault}{\mddefault}{\updefault}{\color[rgb]{0,0,0}$ a^{(1)} $}%
}}}}
\put(1846,-1366){\makebox(0,0)[lb]{\smash{{\SetFigFont{10}{12.0}{\familydefault}{\mddefault}{\updefault}{\color[rgb]{0,0,0}$ a^{(2)} $}%
}}}}
\put(946,-1096){\makebox(0,0)[lb]{\smash{{\SetFigFont{10}{12.0}{\familydefault}{\mddefault}{\updefault}{\color[rgb]{0,0,0}$ a^{(5)} $}%
}}}}
\put(946,-1366){\makebox(0,0)[lb]{\smash{{\SetFigFont{10}{12.0}{\familydefault}{\mddefault}{\updefault}{\color[rgb]{0,0,0}$ a^{(3)} $}%
}}}}
\put(496,-916){\makebox(0,0)[lb]{\smash{{\SetFigFont{10}{12.0}{\familydefault}{\mddefault}{\updefault}{\color[rgb]{0,0,0}$ a^{(4)} $}%
}}}}
\put(3781,-1096){\makebox(0,0)[lb]{\smash{{\SetFigFont{10}{12.0}{\familydefault}{\mddefault}{\updefault}{\color[rgb]{0,0,0}$ V $}%
}}}}
\put(4366,-601){\makebox(0,0)[lb]{\smash{{\SetFigFont{10}{12.0}{\familydefault}{\mddefault}{\updefault}{\color[rgb]{0,0,0}$\Gamma$}%
}}}}
\put(3781,-1591){\makebox(0,0)[lb]{\smash{{\SetFigFont{10}{12.0}{\familydefault}{\mddefault}{\updefault}{\color[rgb]{0,0,0}$2$}%
}}}}
\put(406,-1546){\makebox(0,0)[rb]{\smash{{\SetFigFont{10}{12.0}{\familydefault}{\mddefault}{\updefault}{\color[rgb]{0,0,0}$0$}%
}}}}
\end{picture}%

%% file: pics/lines.tex
\begin{picture}(0,0)%
\includegraphics{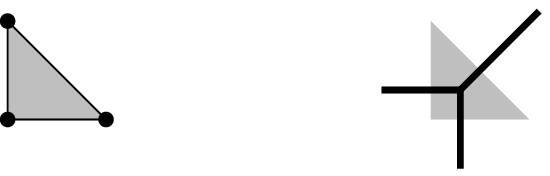}%
\end{picture}%
\setlength{\unitlength}{4144sp}%
\begingroup\makeatletter\ifx\SetFigFont\undefined%
\gdef\SetFigFont#1#2#3#4#5{%
  \reset@font\fontsize{#1}{#2pt}%
  \fontfamily{#3}\fontseries{#4}\fontshape{#5}%
  \selectfont}%
\fi\endgroup%
\begin{picture}(2498,786)(866,-1669)
\put(1081,-1096){\makebox(0,0)[lb]{\smash{{\SetFigFont{10}{12.0}{\familydefault}{\mddefault}{\updefault}{\color[rgb]{0,0,0}$ \Delta_1 $}%
}}}}
\put(3331,-1096){\makebox(0,0)[lb]{\smash{{\SetFigFont{10}{12.0}{\familydefault}{\mddefault}{\updefault}{\color[rgb]{0,0,0}$ \Gamma $}%
}}}}
\end{picture}%

%% file: pics/conics.tex
\begin{picture}(0,0)%
\includegraphics{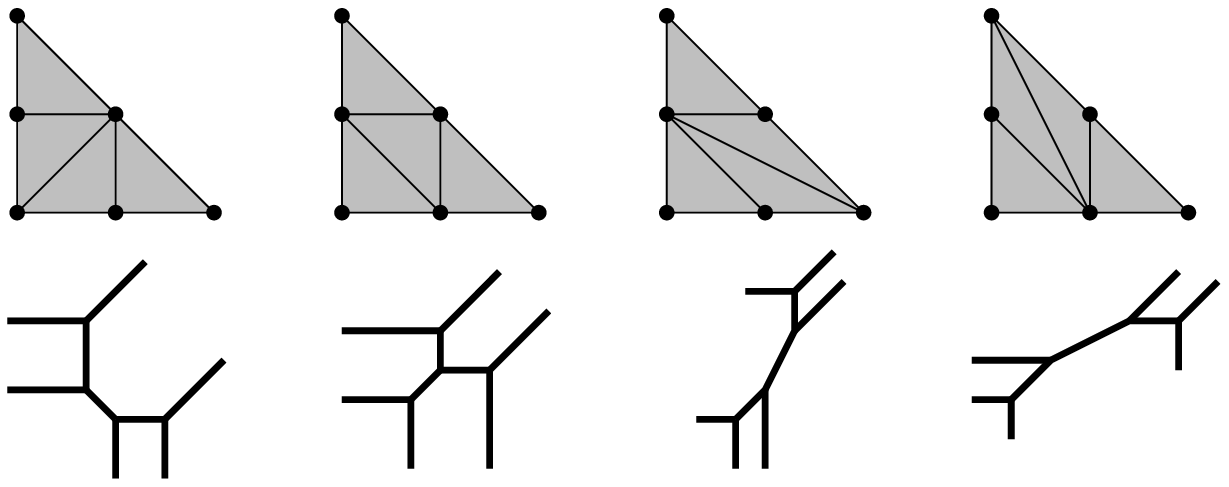}%
\end{picture}%
\setlength{\unitlength}{4144sp}%
\begingroup\makeatletter\ifx\SetFigFont\undefined%
\gdef\SetFigFont#1#2#3#4#5{%
  \reset@font\fontsize{#1}{#2pt}%
  \fontfamily{#3}\fontseries{#4}\fontshape{#5}%
  \selectfont}%
\fi\endgroup%
\begin{picture}(5601,2340)(373,-2366)
\put(991,-1456){\makebox(0,0)[lb]{\smash{{\SetFigFont{10}{12.0}{\familydefault}{\mddefault}{\updefault}{\color[rgb]{0,0,0}$ \Gamma $}%
}}}}
\put(901,-2311){\makebox(0,0)[b]{\smash{{\SetFigFont{10}{12.0}{\familydefault}{\mddefault}{\updefault}{\color[rgb]{0,0,0}(a)}%
}}}}
\put(2386,-2311){\makebox(0,0)[b]{\smash{{\SetFigFont{10}{12.0}{\familydefault}{\mddefault}{\updefault}{\color[rgb]{0,0,0}(b)}%
}}}}
\put(3871,-2311){\makebox(0,0)[b]{\smash{{\SetFigFont{10}{12.0}{\familydefault}{\mddefault}{\updefault}{\color[rgb]{0,0,0}(c)}%
}}}}
\put(5356,-2311){\makebox(0,0)[b]{\smash{{\SetFigFont{10}{12.0}{\familydefault}{\mddefault}{\updefault}{\color[rgb]{0,0,0}(d)}%
}}}}
\put(676,-241){\makebox(0,0)[lb]{\smash{{\SetFigFont{10}{12.0}{\familydefault}{\mddefault}{\updefault}{\color[rgb]{0,0,0}$ \Delta_2 $}%
}}}}
\end{picture}%

%% file: pics/convex.tex
\begin{picture}(0,0)%
\includegraphics{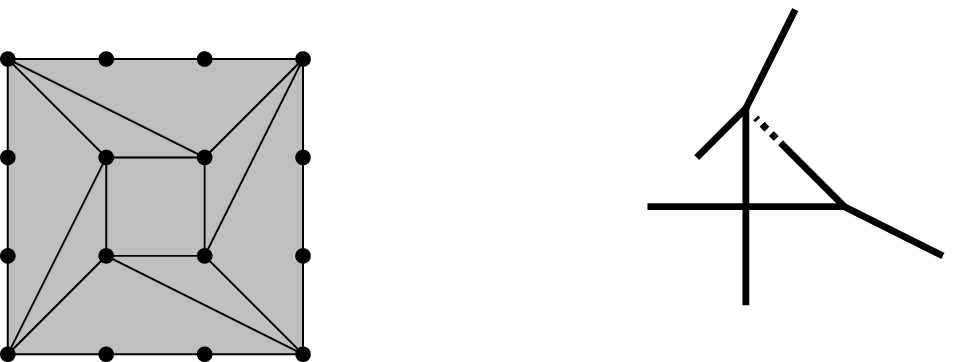}%
\end{picture}%
\setlength{\unitlength}{4144sp}%
\begingroup\makeatletter\ifx\SetFigFont\undefined%
\gdef\SetFigFont#1#2#3#4#5{%
  \reset@font\fontsize{#1}{#2pt}%
  \fontfamily{#3}\fontseries{#4}\fontshape{#5}%
  \selectfont}%
\fi\endgroup%
\begin{picture}(4343,1643)(416,-996)
\put(1846,254){\makebox(0,0)[lb]{\smash{{\SetFigFont{10}{12.0}{\familydefault}{\mddefault}{\updefault}{\color[rgb]{0,0,0}$ \Delta $}%
}}}}
\put(4113,487){\makebox(0,0)[lb]{\smash{{\SetFigFont{10}{12.0}{\familydefault}{\mddefault}{\updefault}{\color[rgb]{0,0,0}$ \Gamma $}%
}}}}
\put(4034,-459){\makebox(0,0)[b]{\smash{{\SetFigFont{10}{12.0}{\familydefault}{\mddefault}{\updefault}{\color[rgb]{0,0,0}$ E_1 $}%
}}}}
\put(3794,-213){\makebox(0,0)[rb]{\smash{{\SetFigFont{10}{12.0}{\familydefault}{\mddefault}{\updefault}{\color[rgb]{0,0,0}$ E_2 $}%
}}}}
\put(3873,347){\makebox(0,0)[rb]{\smash{{\SetFigFont{10}{12.0}{\familydefault}{\mddefault}{\updefault}{\color[rgb]{0,0,0}$ E_3 $}%
}}}}
\put(4127,-99){\makebox(0,0)[lb]{\smash{{\SetFigFont{10}{12.0}{\familydefault}{\mddefault}{\updefault}{\color[rgb]{0,0,0}$ E_4 $}%
}}}}
\put(1127,-200){\makebox(0,0)[b]{\smash{{\SetFigFont{10}{12.0}{\familydefault}{\mddefault}{\updefault}{\color[rgb]{0,0,0}$ E_2 $}%
}}}}
\put(1327,-406){\makebox(0,0)[rb]{\smash{{\SetFigFont{10}{12.0}{\familydefault}{\mddefault}{\updefault}{\color[rgb]{0,0,0}$ E_1 $}%
}}}}
\put(1580,188){\makebox(0,0)[rb]{\smash{{\SetFigFont{10}{12.0}{\familydefault}{\mddefault}{\updefault}{\color[rgb]{0,0,0}$ E_4 $}%
}}}}
\put(1040,127){\makebox(0,0)[lb]{\smash{{\SetFigFont{10}{12.0}{\familydefault}{\mddefault}{\updefault}{\color[rgb]{0,0,0}$ E_3 $}%
}}}}
\end{picture}%

%% file: classical.tex
\section {Tropical versions of classical theorems} \label {classical}


As tropical curves are simply images of classical curves by definition A we can
hope to find tropical --- and thus combinatorial --- versions of many results
known from classical geometry. We will list a few important and interesting
ones in this chapter.

\subsection {Tropical factorization} \label {factor}

Let us start with a very simple statement: in classical geometry it is obvious
that for two polynomials $ f_1,f_2 \in K[z_1,z_2] $ the plane curve defined by
the equation $ (f_1 \cdot f_2)(z) = 0 $ is simply the union of the two curves
with the equations $ f_1(z)=0 $ and $ f_2(z)=0 $.

It is easy to see that an analogous statement holds in the tropical semiring as
well: if $ g_1,g_2 $ are two tropical polynomials, in particular convex
piecewise linear functions, then the corner locus of $ g_1 \odot g_2 = g_1+g_2
$ is simply the union of the corner loci of $ g_1 $ and $ g_2 $. In particular,
the union of two plane tropical curves of degrees $ d_1 $ and $ d_2 $ is always
a plane tropical curve of degree $ d_1+d_2 $.

We can therefore consider the ``tropical factorization problem'' both in a
geometric and an algebraic version. In the geometric language (i.e.\ using
definition C) we would start with a weighted balanced graph in the plane and
ask whether this graph is a union of two weighted subgraphs that are themselves
balanced. In the algebraic language (i.e.\ using definition B) we would start
with a tropical polynomial and ask whether it can be written as a (tropical)
product of two polynomials of smaller degrees.

Note however that these two problems are not entirely equivalent since we have
seen already in section \ref {graphs} that there is no one-to-one
correspondence between tropical polynomials and tropical curves. As an easy
example consider the tropical polynomial $ g(x_1,x_2) = x_1 \oplus x_2 \oplus 0
$ whose corner locus is the curve in figure \ref {pic-limits} (a). If we now
consider the tropical square of this polynomial
\begin {align*}
  g(x_1,x_2) \odot g(x_1,x_2) &= x_1^{\odot 2} \oplus x_2^{\odot 2} \oplus 0
       \oplus (x_1 \odot x_2) \oplus x_1 \oplus x_2 \\
     &= \max \{ 2x_1,2x_2,0,x_1+x_2,x_1,x_2 \}
\end {align*}
then the tropical curve determined by this polynomial is still the same
as before (but with weight 2). But as piecewise linear maps the function $
g(x_1,x_2) \odot g(x_1,x_2) $ is the same as
  \[ \max \{ 2x_1,2x_2,0 \} = x_1^{\odot 2} \oplus x_2^{\odot 2} \oplus 0, \]
and this tropical polynomial \emph {cannot} be written as a product of two
linear tropical polynomials (to prove this just note that there are no additive
inverses in the tropical semiring, so no additive cancellations are possible
anywhere in the expansion of the product).

Geometrically it is in principle easy to decide (although maybe complicated
combinatorially if the degree of the curve is large) whether a given balanced
graph is the union of two smaller ones. Algebraically, it has been shown for
tropical polynomials in one variable that any such polynomial can be replaced
by another one defining the same piecewise linear function that can then be
written as a tropical product of linear factors (see \cite {SS} section 2). For
polynomials in more than one variable not much is known however. The
factorization of a tropical polynomial into irreducible polynomials is in
general not unique, and there is no algorithm known to determine whether a
given polynomial is irreducible resp.\ to compute a (or all) possible
decomposition into irreducible factors (see \cite {SS} section 2). Results in
this direction would be very interesting since it has been shown that a
solution to the tropical factorization problem would also be useful to compute
factorizations of ordinary polynomials more efficiently \cite {GL}.

\subsection {The degree-genus formula} \label {genus}

If $C$ is a smooth complex plane projective curve of degree $d$ then it is
well-known that its genus (i.e.\ the ``number of holes'' in the real surface
$C$) is given by the so-called degree-genus formula $ g = \frac 12 (d-1)(d-2)
$. If $C$ is not smooth then there are several slightly different ways to
define its genus, but for any of these definitions the genus will be at most
the above number $ \frac 12(d-1)(d-2) $.

Let us study the same questions in tropical geometry. If $ \Gamma \subset \RR^2
$ is a plane tropical curve then the most natural way to define its genus is
simply to let it be the number of loops in the graph $ \Gamma $, i.e.\ its
first Betti number $ g= \dim H_1 (\Gamma,\RR) $.

How is this genus related to the degree of $ \Gamma $? To see this let us
denote the set of vertices and bounded edges of $ \Gamma $ by $ \Gamma_0 $ and
$ \Gamma_1 $, respectively. Moreover, for a vertex $ V \in \Gamma_0 $ we define
its \emph {valence} $ \val V $ to be the number of (bounded or unbounded) edges
adjacent to $V$. As $ \Gamma $ has $ 3d $ unbounded and $ |\Gamma_1| $ bounded
edges it follows that
  \[ 3d+2|\Gamma_1| = \sum_{V \in \Gamma_0} \val V. \]
Since the genus of $ \Gamma $ can be computed as $ 1+|\Gamma_1|-|\Gamma_0| $ we
conclude that
  \[ g = 1+\frac 12 \sum_{V \in \Gamma_0} \val V - \frac 32 d - |\Gamma_0|
       = \frac 12 (d-1)(d-2) - \underbrace {\left(
         \frac 12 d^2-\sum_{V \in \Gamma_0} \frac 12(\val V-2) \right)}_{(*)}.
     \]
Recall from section \ref {graphs} that every vertex $ V \in \Gamma_0 $
corresponds to a convex polygon with $ \val V $ vertices in the Newton
subdivision of $ \Delta_d $ corresponding to $ \Gamma $. As such a polygon has
area at least $ \frac 12 (\val V-2) $ and the total area of $ \Delta_d $ is $
\frac 12 d^2 $ it follows that the expression $ (*) $ above is always
non-negative and thus the genus of $ \Gamma $ is always at most $ \frac 12
(d-1)(d-2) $, as in the classical case. Equality holds if and only if all
polygons in the Newton subdivision have minimal area for its number of
vertices.

There is nothing like a general ``tropical singularity theory'' yet, but
usually one says that a plane tropical curve $ \Gamma $ is \emph {smooth} if
its Newton subdivision is maximal (i.e.\ consists of $ d^2 $ triangles of area
$ \frac 12 $ each), or equivalently if every vertex of $ \Gamma $ has valence
3, all weights of the edges are 1, and the primitive integral vectors along the
edges adjacent to any vertex generate the lattice $ \ZZ^2 $. With this
definition it follows from our computations above that \emph {the genus of a
smooth plane tropical curve of degree $d$ is $ \frac 12 (d-1)(d-2) $}, just as
in classical geometry.

The following picture shows three examples of plane cubics: the curve (a) is
smooth (and hence of genus 1), (b) is not smooth but still of genus 1 (since
the Newton subdivision contains only a parallelogram of area $1$ and triangles
of area $ \frac 12 $), and (c) has genus 0 (since the Newton subdivision
contains triangles of area greater than $ \frac 12 $).

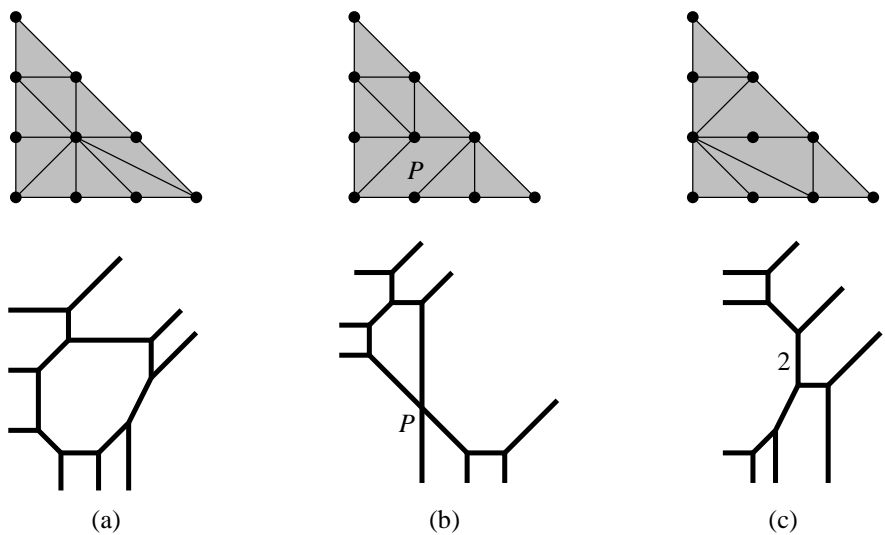
\begin {figure}[H]
  \begin {center} \input {pics/cubics} \end {center}
  \caption {Three plane tropical cubics}
  \label {pic-cubics}
\end {figure}

Let us consider case (b) in more detail. Obviously the parallelogram $P$ in the
Newton subdivision gives rise to a point in the tropical curve (that we also
denoted by $P$) where two straight edges intersect. We can therefore think of $
\Gamma $ as the planar image of a graph of genus 0 that has a ``crossing'' at
$P$. This corresponds exactly to a normal crossing singularity in the classical
case, i.e.\ to a complex curve $C$ with a point $ P \in C $ where two smooth
branches meet transversely. In fact, for a plane tropical curve whose Newton
subdivision contains only of triangles and parallelograms one sometimes
subtracts the number of parallelograms from the genus defined above for that
reason (so that e.g.\ the genus of the curve (b) above would then be 0).

\subsection {B\'ezout's theorem} \label {bezout}

Let $ C_1 $ and $ C_2 $ be two distinct smooth plane projective complex curves
of degrees $ d_1 $ and $ d_2 $, respectively. B\'ezout's theorem states that
the intersection $ C_1 \cap C_2 $ then consists of at most $ d_1 \cdot d_2 $
points, and that in fact equality holds if the intersection points are counted
with the correct multiplicity, namely with the local intersection
multiplicities of $ C_1 $ and $ C_2 $.

We have seen in section \ref {generalizations} already that there is a slight
problem if we try to find an analogous statement in tropical geometry: it is
not even true that the intersection of two distinct (smooth) plane tropical
curves is always finite --- they might as well share some common line segments.

Let us ignore this problem for a moment however and assume that we have two
smooth tropical curves $ \Gamma_1 $ and $ \Gamma_2 $ of degrees $ d_1 $ and $
d_2 $ respectively that intersect in finitely many points, and that none of
these intersection points is a vertex of either curve. In this case it is in
fact easy to find a tropical B\'ezout theorem: by section \ref {factor} the
union $ \Gamma_1 \cup \Gamma_2 $ is a plane tropical curve of degree $ d_1+d_2
$ and hence corresponds to a Newton subdivision of $ \Delta_{d_1+d_2} $. The
vertices of $ \Gamma_1 \cup \Gamma_2 $ are of two types:
\begin {itemize}
\item the 3-valent vertices of $ \Gamma_1 $ and $ \Gamma_2 $ are of course also
  present in $ \Gamma_1 \cup \Gamma_2 $. As $ \Gamma_1 \cup \Gamma_2 $ looks
  locally the same as $ \Gamma_1 $ resp.\ $ \Gamma_2 $ around such a vertex the
  triangles in the Newton subdivisions for $ \Gamma_1 $ and $ \Gamma_2 $
  can also be found in the subdivision for $ \Gamma_1 \cup \Gamma_2 $.
\item every intersection point in $ \Gamma_1 \cap \Gamma_2 $ gives rise to a
  4-valent vertex of $ \Gamma_1 \cup \Gamma_2 $ where two straight lines meet.
  As in the end of section \ref {genus} this gives rise to a parallelogram in
  the Newton subdivision of $ \Gamma_1 \cup \Gamma_2 $.
\end {itemize}
The following picture illustrates this for the case of two conics, one of type
(a) and one of type (d) in the notation of figure \ref {pic-conics}:

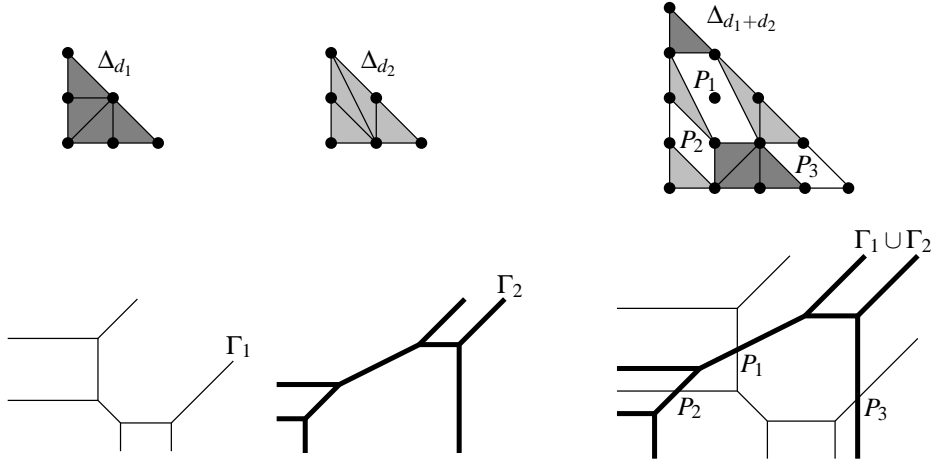
\begin {figure}[H]
  \begin {center} \input {pics/bezout} \end {center}
  \caption {The intersection of two smooth plane conics}
  \label {pic-bezout}
\end {figure}

We conclude that the area covered by the parallelograms corresponding to the
points in $ \Gamma_1 \cap \Gamma_2 $ is
  \[ \area (\Delta_{d_1+d_2}) - \area \Delta_{d_1} - \area \Delta_{d_2}
       = \frac 12 (d_1+d_2)^2 - \frac 12 d_1^2 - \frac 12 d_2^2 = d_1d_2. \]
So if we define the intersection multiplicity of $ \Gamma_1 $ and $ \Gamma_2 $
in a point $ P \in \Gamma_1 \cap \Gamma_2 $ to be the area of the parallelogram
corresponding to $P$ in the Newton subdivision of $ \Delta_{d_1+d_2} $ (which
is a positive integer) then it follows that the sum of all these intersection
multiplicities is $ d_1d_2 $. This is the tropical B\'ezout theorem that has
first appeared in the literature in \cite {RST} (albeit with a different
proof). For example, in figure \ref {pic-bezout} there are three intersection
points of which $ P_1 $ has multiplicity 2 and the other two have multiplicity
1, so the total weighted number is $ 4 = \deg \Gamma_1 \cdot \deg \Gamma_2 $.

It is easy to see that one can also interpret the intersection multiplicity in
terms of the graphs $ \Gamma_1 $ and $ \Gamma_2 $ themselves, without using
their Newton subdivisions: if $ w_1, w_2 $ are the weights and $ u_1,u_2 $ the
primitive integral vectors along the two edges meeting in a point $ P \in
\Gamma_1 \cap \Gamma_2 $ then the intersection multiplicity of $ \Gamma_1 $ and
$ \Gamma_2 $ in $P$ is $ w_1 \cdot w_2 \cdot |\det (u_1,u_2)| $.

It is quite straightforward to generalize our tropical B\'ezout theorem to
other cases. For example, it is easily checked that our genericity assumptions
(i.e.\ that the curves are smooth and do not intersect in vertices) are not
necessary --- with a suitable definition of intersection multiplicity the same
theorem holds if only the intersection $ \Gamma_1 \cap \Gamma_2 $ is finite.
Moreover, by replacing the triangles $ \Delta_d $ by other convex polytopes the
very same proof can be used for tropical curves in toric surfaces (see section
\ref {generalizations}), leading to some very basic intersection theory for
curves on such surfaces. The above proof can also be iterated to the case of
hypersurfaces: if $n$ tropical hypersurfaces in $ \RR^n $ intersect in finitely
many points then the weighted number of such points is the product of the
degrees of the hypersurfaces (resp.\ the mixed volume of the convex polytopes
in the case of hypersurfaces in toric varieties). There is however no tropical
version yet of Chow groups and a general intersection theory in the sense of
Fulton \cite {F}.

Much more surprising is the fact that --- unlike in classical geometry --- one
can in fact get a version of B\'ezout's theorem without \emph {any} condition
on the curves $ \Gamma_1 $ and $ \Gamma_2 $, even if they coincide or share
common line segments. In this case the strategy is to move one of the curves,
say $ \Gamma_2 $, to a nearby curve $ \Gamma_2' $ so that the intersection is
finite. As $ \Gamma_2' $ is moved back to $ \Gamma_2 $ it can be shown that the
finitely many intersection points in $ \Gamma_1 \cap \Gamma_2' $ have
well-defined limit points in $ \Gamma_1 \cap \Gamma_2 $ (see \cite {RST}
section 4). This is called the \emph {stable intersection} of $ \Gamma_1 $ and
$ \Gamma_2 $. For example, in the following picture the stable intersection of
the conic $ \Gamma_1 $ with itself is simply the union of the four vertices $
Q_1,\dots,Q_4 $ of the curve, each counted with multiplicity 1.

\begin {figure}[H]
  \begin {center} \input {pics/stable} \end {center}
  \caption {The stable intersection of a conic with itself}
  \label {pic-stable}
\end {figure}
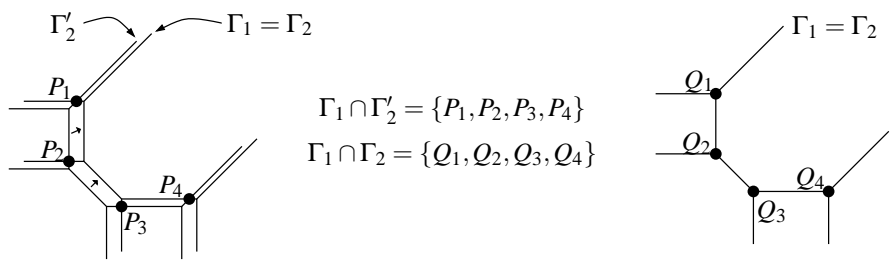

The proof of B\'ezout's theorem that we have given above is purely
combinatorial and does not use the corresponding classical statement. In fact,
we have mentioned already that one of the big advantages of tropical geometry
is that possibly complicated algebraic or geometric questions can be reduced to
entirely combinatorial ones. Sometimes it is also interesting and instructive
however to remember that tropical geometry is nothing but an ``image of
classical geometry (over the field $K$ of Puiseux series) under the valuation
map''. This way one can try to transfer results from classical geometry
directly to tropical geometry. In our case at hand we could proceed as follows:
if $ \Gamma_1 $ and $ \Gamma_2 $ are two plane tropical curves of degrees $ d_1
$ and $ d_2 $ respectively they can be realized as the images under the
valuation map of two classical plane curves $ C_1 $ and $ C_2 $ of these
degrees over $K$. Now if we use the classical B\'ezout theorem for these curves
we can conclude that $ C_1 $ and $ C_2 $ intersect in $ d_1 d_2 $ points
(counted with the correct multiplicities). Of course the images of these points
under the valuation map are intersection points of $ \Gamma_1 \cap \Gamma_2 $,
and in a generic situation these will be the only intersection points of these
tropical curves. To make B\'ezout's theorem hold in the tropical setting we
therefore should define the intersection multiplicity of $ \Gamma_1 $ and $
\Gamma_2 $ in a point $P$ to be the sum of the intersection multiplicities of $
C_1 $ and $ C_2 $ at all points that map to $P$ under the valuation map.

Let us check in a simple example that this agrees with the definition of
tropical intersection multiplicity that we have given above. We consider a
local situation around an intersection point of $ \Gamma_1 $ and $ \Gamma_2 $.
For simplicity let us choose coordinates so that the intersection point is the
origin and the two curves have local equations $ \Gamma_1 = \{ (x_1,x_2) ;\;
x_2 = 0 \} $ and $ \Gamma_2 = \{ (x_1,x_2);\; x_2 = n x_1 \} $ for some $ n \in
\NN_{>0} $. We can then choose plane curves $ C_1 $ and $ C_2 $ mapping to $
\Gamma_1 $ and $ \Gamma_2 $ as in the following picture:

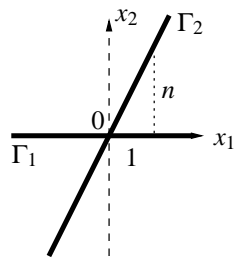
\begin {figure}[H]
  \begin {center} \input {pics/puiseux} \end {center}
  \caption {Classical interpretation of the tropical intersection multiplicity}
  \label {pic-puiseux}
\end {figure}

We see that the intersection $ C_1 \cap C_2 $ consists of $n$ points ---
corresponding to a choice of $n$-th root of unity for $ z_1 $ --- that all map
to $ (0,0) $ in the tropical picture under the valuation map. So the
intersection multiplicity of $ \Gamma_1 $ and $ \Gamma_2 $ should be defined to
be $n$. In fact this agrees with our definition above since the primitive
integral vectors of the two tropical curves are $ u_1 = (1,0) $ and $ u_2 (1,n)
$, and $ |\det (u_1,u_2)| = n $.

\subsection {The group structure of a plane cubic curve} \label {elliptic}

Again let $C$ be a smooth plane projective complex curve. The group of \emph
{divisors} $ \Div C $ on $C$ is the free abelian group generated by the points
of $C$, i.e.\ the points of $ \Div C $ are finite formal linear combinations $
D = \lambda_1 P_1 + \cdots + \lambda_n P_n $ with $ \lambda_i \in \ZZ $ and $
P_i \in C $. For such a divisor we call the number $ \lambda_1 + \cdots +
\lambda_n \in \ZZ $ the degree of $D$. Obviously the subset $ \Div^0(C) $ of
divisors of degree 0 is a subgroup of $ \Div C $.

If $ C' $ is any other curve then the intersection $ C \cap C' $ consists of
finitely many points. Hence we can consider this intersection to be an element
of $ \Div C $, where we count each point with its intersection multiplicity.
This divisor is usually denoted $ C \cdot C' $. By B\'ezout's theorem its
degree is $ \deg C \cdot \deg C' $.

Two divisors $ D_1,D_2 \in \Div C $ are called \emph {equivalent} if there
are curves $ C',C'' $ of the same degree such that $ D_1 - D_2 = C \cdot C' - C
\cdot C'' $. The group of equivalence classes is usually denoted $ \Pic (C) $,
or $ \Pic^0 (C) $ if we restrict to divisors of degree 0.

Now if $C$ has degree 3 it can be shown that after picking a base point $ P_0
\in C $ the map $ C \to \Pic^0(C), P \mapsto P-P_0 $ is a bijection. Hence we
can use this map to define a group structure on $C$, or vice versa the
structure of a complex algebraic curve on $ \Pic^0(C) $. Alternatively, it
follows by the degree-genus formula that $C$ has genus 1, i.e.\ it is a torus.
It can be shown that one can realize this torus in the form $ \CC/\Lambda $
with a lattice $ \Lambda \cong \ZZ^2 $ so that the group structure on $
\CC/\Lambda $ induced by addition on $ \CC $ is precisely the group structure
defined by the bijection with $ \Pic^0 (C) $ constructed above.

Which of all these results remain true in the tropical world? Let $ \Gamma $ be
a plane tropical curve, and let us start by defining the group of divisors $
\Div \Gamma $ in the same way as above, i.e.\ as the free abelian group
generated by the points of $ \Gamma $. As we have seen already that B\'ezout's
theorem holds in the tropical set-up as well we can also define the group $
\Pic^0 (\Gamma) $ in the same way as before, i.e.\ it is the group of divisors
on $ \Gamma $ of degree 0 modulo those that can be written as $ \Gamma \cdot
\Gamma' - \Gamma \cdot \Gamma'' $ for some tropical curves $ \Gamma' $ and $
\Gamma'' $ of the same degree.

We can also still define a map $ \Gamma \to \Pic^0 (\Gamma) $ by sending $P$ to
$ P-P_0 $ after choosing a base point $ P_0 \in P $. However, this map is \emph
{not} a bijection in the tropical case as can be seen from the following
picture:

\begin {figure}[H]
  \begin {center} \input {pics/group} \end {center}
  \caption {The group $ \Pic^0 (\Gamma) $ of a plane cubic curve $ \Gamma $}
  \label {pic-group}
\end {figure}
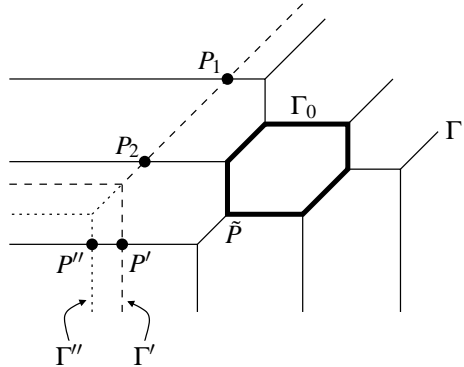

As the two lines $ \Gamma' $ and $ \Gamma'' $ meet the cubic curve $ \Gamma $
in the points $ P_1,P_2,P' $ and $ P_1,P_2,P'' $, respectively, it follows that
$ P' $ and $ P'' $ are equivalent; in particular $ P' - P_0 = P'' - P_0 \in
\Pic^0(\Gamma) $ for any choice of base point $ P_0 $.

In fact, if $ \Gamma_0 $ is the unique loop in $ \Gamma $ (drawn in bold in the
picture above) Vigeland \cite {V} has shown by the same techniques that any
point of $ \Gamma \backslash \Gamma_0 $ is equivalent to the point of $
\Gamma_0 $ ``nearest to it'', e.g.\ $ P' $ and $ P'' $ in figure \ref
{pic-group} are both equivalent to $ \tilde P $. Moreover, he has proven that
the map $ P \mapsto P-P_0 $ in fact gives a bijection of $ \Gamma_0 $ with $
\Pic^0(\Gamma) $, so that the loop $ \Gamma_0 $ inherits a group structure from
$ \Pic^0(\Gamma) $. In fact, this group structure is just the ordinary group
structure of the unit circle $ S^1 $ after choosing a suitable bijection of $
S^1 $ with $ \Gamma_0 $.

So to a certain extent we can find analogues of the classical results about
plane cubics in the tropical world as well. Despite these encouraging results
one should note however that there is no good theory of divisors and their
equivalence yet on arbitrary tropical curves.

%% file: pics/cubics.tex
\begin{picture}(0,0)%
\includegraphics{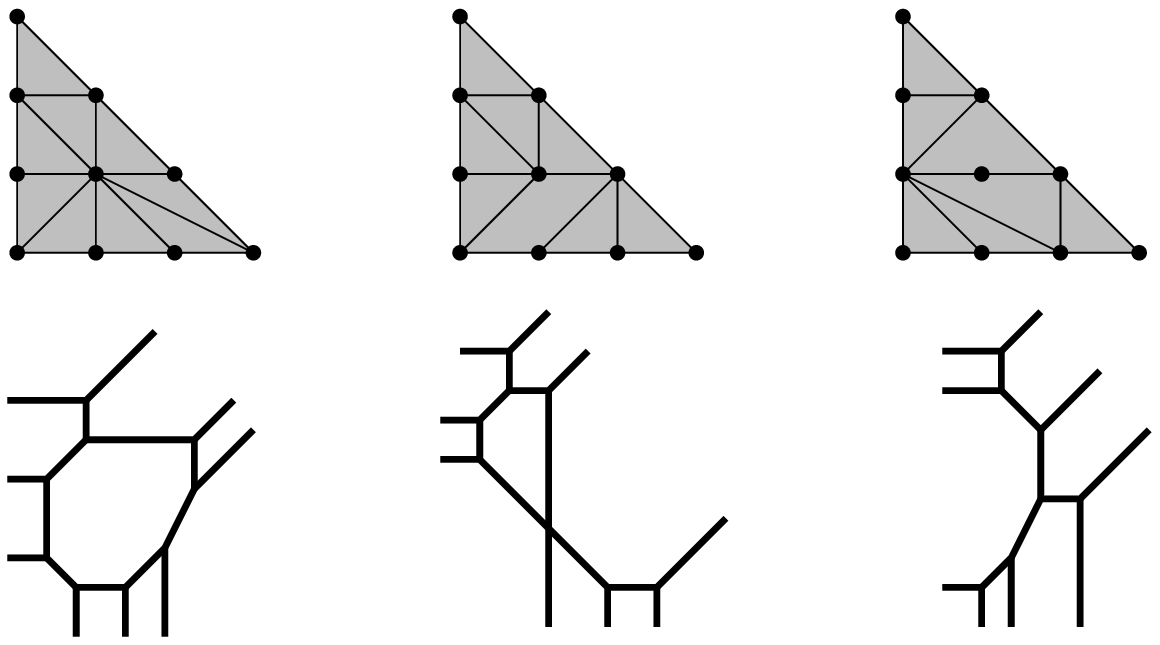}%
\end{picture}%
\setlength{\unitlength}{4144sp}%
\begingroup\makeatletter\ifx\SetFigFont\undefined%
\gdef\SetFigFont#1#2#3#4#5{%
  \reset@font\fontsize{#1}{#2pt}%
  \fontfamily{#3}\fontseries{#4}\fontshape{#5}%
  \selectfont}%
\fi\endgroup%
\begin{picture}(5286,3150)(373,-2726)
\put(5041,-2671){\makebox(0,0)[b]{\smash{{\SetFigFont{10}{12.0}{\familydefault}{\mddefault}{\updefault}{\color[rgb]{0,0,0}(c)}%
}}}}
\put(3016,-2671){\makebox(0,0)[b]{\smash{{\SetFigFont{10}{12.0}{\familydefault}{\mddefault}{\updefault}{\color[rgb]{0,0,0}(b)}%
}}}}
\put(991,-2671){\makebox(0,0)[b]{\smash{{\SetFigFont{10}{12.0}{\familydefault}{\mddefault}{\updefault}{\color[rgb]{0,0,0}(a)}%
}}}}
\put(5086,-1726){\makebox(0,0)[rb]{\smash{{\SetFigFont{10}{12.0}{\familydefault}{\mddefault}{\updefault}{\color[rgb]{0,0,0}$2$}%
}}}}
\put(2842,-583){\makebox(0,0)[b]{\smash{{\SetFigFont{10}{12.0}{\familydefault}{\mddefault}{\updefault}{\color[rgb]{0,0,0}$P$}%
}}}}
\put(2838,-2096){\makebox(0,0)[rb]{\smash{{\SetFigFont{10}{12.0}{\familydefault}{\mddefault}{\updefault}{\color[rgb]{0,0,0}$P$}%
}}}}
\end{picture}%

%% file: pics/bezout.tex
\begin{picture}(0,0)%
\includegraphics{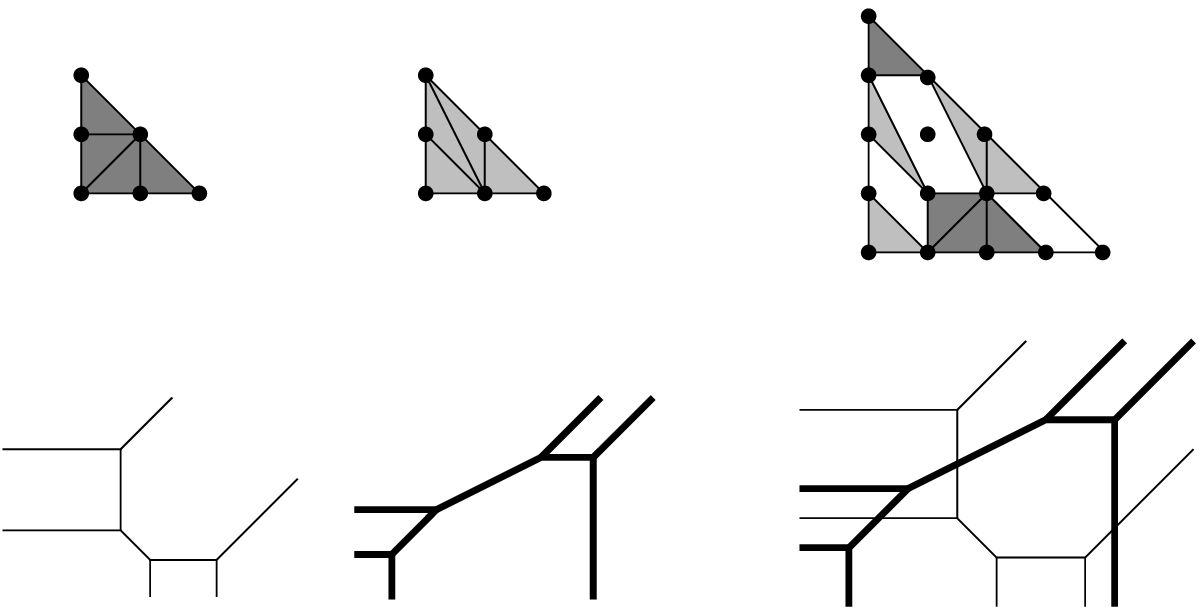}%
\end{picture}%
\setlength{\unitlength}{4144sp}%
\begingroup\makeatletter\ifx\SetFigFont\undefined%
\gdef\SetFigFont#1#2#3#4#5{%
  \reset@font\fontsize{#1}{#2pt}%
  \fontfamily{#3}\fontseries{#4}\fontshape{#5}%
  \selectfont}%
\fi\endgroup%
\begin{picture}(5562,2790)(-236,-2074)
\put(3961,569){\makebox(0,0)[lb]{\smash{{\SetFigFont{10}{12.0}{\familydefault}{\mddefault}{\updefault}{\color[rgb]{0,0,0}$ \Delta_{d_1+d_2} $}%
}}}}
\put(4486,-334){\makebox(0,0)[lb]{\smash{{\SetFigFont{10}{12.0}{\familydefault}{\mddefault}{\updefault}{\color[rgb]{0,0,0}$ P_3 $}%
}}}}
\put(3878,181){\makebox(0,0)[lb]{\smash{{\SetFigFont{10}{12.0}{\familydefault}{\mddefault}{\updefault}{\color[rgb]{0,0,0}$ P_1 $}%
}}}}
\put(3797,-183){\makebox(0,0)[lb]{\smash{{\SetFigFont{10}{12.0}{\familydefault}{\mddefault}{\updefault}{\color[rgb]{0,0,0}$ P_2 $}%
}}}}
\put(5311,-781){\makebox(0,0)[rb]{\smash{{\SetFigFont{10}{12.0}{\familydefault}{\mddefault}{\updefault}{\color[rgb]{0,0,0}$ \Gamma_1 \cup \Gamma_2 $}%
}}}}
\put(4158,-1524){\makebox(0,0)[lb]{\smash{{\SetFigFont{10}{12.0}{\familydefault}{\mddefault}{\updefault}{\color[rgb]{0,0,0}$ P_1 $}%
}}}}
\put(3783,-1777){\makebox(0,0)[lb]{\smash{{\SetFigFont{10}{12.0}{\familydefault}{\mddefault}{\updefault}{\color[rgb]{0,0,0}$ P_2 $}%
}}}}
\put(4899,-1787){\makebox(0,0)[lb]{\smash{{\SetFigFont{10}{12.0}{\familydefault}{\mddefault}{\updefault}{\color[rgb]{0,0,0}$ P_3 $}%
}}}}
\put(316,299){\makebox(0,0)[lb]{\smash{{\SetFigFont{10}{12.0}{\familydefault}{\mddefault}{\updefault}{\color[rgb]{0,0,0}$ \Delta_{d_1} $}%
}}}}
\put(1891,299){\makebox(0,0)[lb]{\smash{{\SetFigFont{10}{12.0}{\familydefault}{\mddefault}{\updefault}{\color[rgb]{0,0,0}$ \Delta_{d_2} $}%
}}}}
\put(1081,-1411){\makebox(0,0)[lb]{\smash{{\SetFigFont{10}{12.0}{\familydefault}{\mddefault}{\updefault}{\color[rgb]{0,0,0}$ \Gamma_1 $}%
}}}}
\put(2701,-1051){\makebox(0,0)[lb]{\smash{{\SetFigFont{10}{12.0}{\familydefault}{\mddefault}{\updefault}{\color[rgb]{0,0,0}$ \Gamma_2 $}%
}}}}
\end{picture}%

%% file: pics/stable.tex
\begin{picture}(0,0)%
\includegraphics{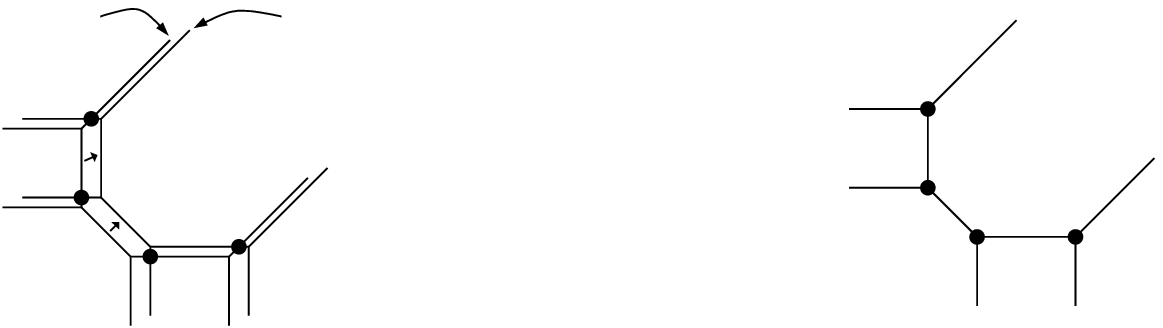}%
\end{picture}%
\setlength{\unitlength}{4144sp}%
\begingroup\makeatletter\ifx\SetFigFont\undefined%
\gdef\SetFigFont#1#2#3#4#5{%
  \reset@font\fontsize{#1}{#2pt}%
  \fontfamily{#3}\fontseries{#4}\fontshape{#5}%
  \selectfont}%
\fi\endgroup%
\begin{picture}(5289,1520)(349,-973)
\put(3016,-376){\makebox(0,0)[b]{\smash{{\SetFigFont{10}{12.0}{\familydefault}{\mddefault}{\updefault}{\color[rgb]{0,0,0}$ \Gamma_1 \cap \Gamma_2 = \{Q_1,Q_2,Q_3,Q_4\} $}%
}}}}
\put(3016,-106){\makebox(0,0)[b]{\smash{{\SetFigFont{10}{12.0}{\familydefault}{\mddefault}{\updefault}{\color[rgb]{0,0,0}$ \Gamma_1 \cap \Gamma_2' = \{P_1,P_2,P_3,P_4\} $}%
}}}}
\put(1667,400){\makebox(0,0)[lb]{\smash{{\SetFigFont{10}{12.0}{\familydefault}{\mddefault}{\updefault}{\color[rgb]{0,0,0}$ \Gamma_1 = \Gamma_2 $}%
}}}}
\put(771,392){\makebox(0,0)[rb]{\smash{{\SetFigFont{10}{12.0}{\familydefault}{\mddefault}{\updefault}{\color[rgb]{0,0,0}$ \Gamma_2' $}%
}}}}
\put(5047,387){\makebox(0,0)[lb]{\smash{{\SetFigFont{10}{12.0}{\familydefault}{\mddefault}{\updefault}{\color[rgb]{0,0,0}$ \Gamma_1 = \Gamma_2 $}%
}}}}
\put(1051,-789){\makebox(0,0)[lb]{\smash{{\SetFigFont{10}{12.0}{\familydefault}{\mddefault}{\updefault}{\color[rgb]{0,0,0}$ P_3 $}%
}}}}
\put(701,-354){\makebox(0,0)[rb]{\smash{{\SetFigFont{10}{12.0}{\familydefault}{\mddefault}{\updefault}{\color[rgb]{0,0,0}$ P_2 $}%
}}}}
\put(745, 19){\makebox(0,0)[rb]{\smash{{\SetFigFont{10}{12.0}{\familydefault}{\mddefault}{\updefault}{\color[rgb]{0,0,0}$ P_1 $}%
}}}}
\put(1410,-567){\makebox(0,0)[rb]{\smash{{\SetFigFont{10}{12.0}{\familydefault}{\mddefault}{\updefault}{\color[rgb]{0,0,0}$ P_4 $}%
}}}}
\put(4585, 66){\makebox(0,0)[rb]{\smash{{\SetFigFont{10}{12.0}{\familydefault}{\mddefault}{\updefault}{\color[rgb]{0,0,0}$ Q_1 $}%
}}}}
\put(4568,-307){\makebox(0,0)[rb]{\smash{{\SetFigFont{10}{12.0}{\familydefault}{\mddefault}{\updefault}{\color[rgb]{0,0,0}$ Q_2 $}%
}}}}
\put(5250,-520){\makebox(0,0)[rb]{\smash{{\SetFigFont{10}{12.0}{\familydefault}{\mddefault}{\updefault}{\color[rgb]{0,0,0}$ Q_4 $}%
}}}}
\put(4838,-715){\makebox(0,0)[lb]{\smash{{\SetFigFont{10}{12.0}{\familydefault}{\mddefault}{\updefault}{\color[rgb]{0,0,0}$ Q_3 $}%
}}}}
\end{picture}%

%% file: pics/puiseux.tex
\begin{picture}(0,0)%
\includegraphics{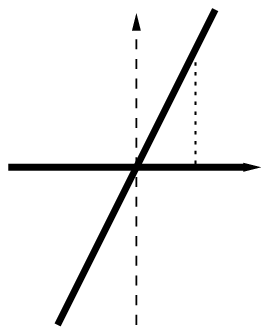}%
\end{picture}%
\setlength{\unitlength}{4144sp}%
\begingroup\makeatletter\ifx\SetFigFont\undefined%
\gdef\SetFigFont#1#2#3#4#5{%
  \reset@font\fontsize{#1}{#2pt}%
  \fontfamily{#3}\fontseries{#4}\fontshape{#5}%
  \selectfont}%
\fi\endgroup%
\begin{picture}(4032,1585)(1336,-904)
\put(4861,-331){\makebox(0,0)[b]{\smash{{\SetFigFont{10}{12.0}{\familydefault}{\mddefault}{\updefault}{\color[rgb]{0,0,0}$1$}%
}}}}
\put(5041, 74){\makebox(0,0)[lb]{\smash{{\SetFigFont{10}{12.0}{\familydefault}{\mddefault}{\updefault}{\color[rgb]{0,0,0}$n$}%
}}}}
\put(5131,479){\makebox(0,0)[lb]{\smash{{\SetFigFont{10}{12.0}{\familydefault}{\mddefault}{\updefault}{\color[rgb]{0,0,0}$ \Gamma_2 $}%
}}}}
\put(5353,-208){\makebox(0,0)[lb]{\smash{{\SetFigFont{10}{12.0}{\familydefault}{\mddefault}{\updefault}{\color[rgb]{0,0,0}$ x_1 $}%
}}}}
\put(4765,534){\makebox(0,0)[lb]{\smash{{\SetFigFont{10}{12.0}{\familydefault}{\mddefault}{\updefault}{\color[rgb]{0,0,0}$ x_2 $}%
}}}}
\put(4698,-113){\makebox(0,0)[rb]{\smash{{\SetFigFont{10}{12.0}{\familydefault}{\mddefault}{\updefault}{\color[rgb]{0,0,0}$0$}%
}}}}
\put(1351,-16){\makebox(0,0)[lb]{\smash{{\SetFigFont{10}{12.0}{\familydefault}{\mddefault}{\updefault}{\color[rgb]{0,0,0}$ \Gamma_2 = \{ (x_1,x_2) \in \RR^2;\; x_2 = nx_1 \} $}%
}}}}
\put(1351,-466){\makebox(0,0)[lb]{\smash{{\SetFigFont{10}{12.0}{\familydefault}{\mddefault}{\updefault}{\color[rgb]{0,0,0}$ C_1 = \{ (z_1,z_2) \in K^2 ;\; z_2 = 1 \} $}%
}}}}
\put(1351,-691){\makebox(0,0)[lb]{\smash{{\SetFigFont{10}{12.0}{\familydefault}{\mddefault}{\updefault}{\color[rgb]{0,0,0}$ C_2 = \{ (z_1,z_2) \in K^2 ;\; z_2 = z_1^n \} $}%
}}}}
\put(1351,209){\makebox(0,0)[lb]{\smash{{\SetFigFont{10}{12.0}{\familydefault}{\mddefault}{\updefault}{\color[rgb]{0,0,0}$ \Gamma_1 = \{ (x_1,x_2) \in \RR^2;\; x_2 = 0 \} $}%
}}}}
\put(4141,-312){\makebox(0,0)[lb]{\smash{{\SetFigFont{10}{12.0}{\familydefault}{\mddefault}{\updefault}{\color[rgb]{0,0,0}$ \Gamma_1 $}%
}}}}
\end{picture}%

%% file: pics/group.tex
\begin{picture}(0,0)%
\includegraphics{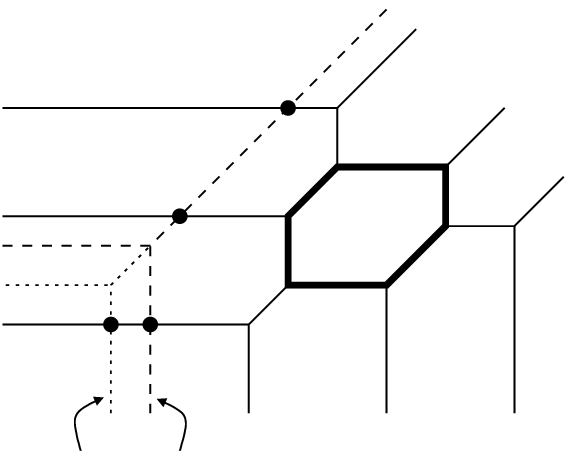}%
\end{picture}%
\setlength{\unitlength}{4144sp}%
\begingroup\makeatletter\ifx\SetFigFont\undefined%
\gdef\SetFigFont#1#2#3#4#5{%
  \reset@font\fontsize{#1}{#2pt}%
  \fontfamily{#3}\fontseries{#4}\fontshape{#5}%
  \selectfont}%
\fi\endgroup%
\begin{picture}(2637,2227)(1114,-1736)
\put(3736,-331){\makebox(0,0)[lb]{\smash{{\SetFigFont{10}{12.0}{\familydefault}{\mddefault}{\updefault}{\color[rgb]{0,0,0}$ \Gamma $}%
}}}}
\put(1936,-1681){\makebox(0,0)[b]{\smash{{\SetFigFont{10}{12.0}{\familydefault}{\mddefault}{\updefault}{\color[rgb]{0,0,0}$ \Gamma' $}%
}}}}
\put(1486,-1681){\makebox(0,0)[b]{\smash{{\SetFigFont{10}{12.0}{\familydefault}{\mddefault}{\updefault}{\color[rgb]{0,0,0}$ \Gamma'' $}%
}}}}
\put(2403, 89){\makebox(0,0)[rb]{\smash{{\SetFigFont{10}{12.0}{\familydefault}{\mddefault}{\updefault}{\color[rgb]{0,0,0}$ P_1 $}%
}}}}
\put(1910,-418){\makebox(0,0)[rb]{\smash{{\SetFigFont{10}{12.0}{\familydefault}{\mddefault}{\updefault}{\color[rgb]{0,0,0}$ P_2 $}%
}}}}
\put(1577,-1130){\makebox(0,0)[rb]{\smash{{\SetFigFont{10}{12.0}{\familydefault}{\mddefault}{\updefault}{\color[rgb]{0,0,0}$ P'' $}%
}}}}
\put(1844,-1125){\makebox(0,0)[lb]{\smash{{\SetFigFont{10}{12.0}{\familydefault}{\mddefault}{\updefault}{\color[rgb]{0,0,0}$ P' $}%
}}}}
\put(2420,-953){\makebox(0,0)[lb]{\smash{{\SetFigFont{10}{12.0}{\familydefault}{\mddefault}{\updefault}{\color[rgb]{0,0,0}$ \tilde P $}%
}}}}
\put(2809,-184){\makebox(0,0)[lb]{\smash{{\SetFigFont{10}{12.0}{\familydefault}{\mddefault}{\updefault}{\color[rgb]{0,0,0}$ \Gamma_0 $}%
}}}}
\end{picture}%

%% file: enum.tex
\section {Tropical techniques in enumerative geometry} \label {enum}


In the last chapter we have seen that many classical results from algebraic
geometry have a tropical counterpart. However, the main reason why tropical
geometry received so much attention recently is that it can be used very
successfully to solve even complicated problems in complex and real enumerative
geometry. So in the rest of this paper we want to give a brief sketch of the
progress that has been made so far in tropical enumerative geometry.

\subsection {Complex enumerative geometry and Gromov-Witten invariants}
  \label {gromov-witten}

If we stick to plane curves the main basic question in complex enumerative
geometry is: given $ d \ge 1 $ and $ g \ge 0 $, what are the numbers $ N_{g,d}
$ of curves of genus $g$ and degree $d$ in the complex projective plane that
pass through $ 3d+g-1 $ general given points? (The number $ 3d+g-1 $ is chosen
so that a naive count of dimensions versus conditions leads one to expect a
finite non-zero answer.) Except for some special cases the answer to this
problem has not been known until the invention of Gromov-Witten theory about
ten years ago.

For curves in projective spaces the main objects of study in Gromov-Witten
theory are the so-called \emph {moduli spaces of stable maps} $ \bar M_{g,n}
(r,d) $ for $ n,r \ge 0 $. We will be more specific about the definition of
these spaces in section \ref {wdvv} --- for the moment it suffices to say that
they are reasonably well-behaved, compact spaces that parametrize curves of
genus $g$ and degree $d$ with $n$ marked points in $ \PP^r $. Of course there
are \emph {evaluation maps} $ \ev_i: \bar M_{g,n} (r,d) \to \PP^r $ for $
i=1,\dots,n $ that map such an $n$-pointed curve to the position of its $i$-th
marked point in $ \PP^r $.

In the case of plane curves mentioned above we now set $ n=3d+g-1 $ and choose
$n$ general points $ P_1,\dots,P_n \in \PP^2 $. The intersection $ \ev_1^{-1}
P_1 \cap \cdots \cap \ev_n^{-1} P_n $ then obviously corresponds to those
plane curves of the given genus and degree that pass through the specified
points. As we have chosen $n$ so that the expected dimension of this
intersection is 0 it makes sense to \emph {define} the number $ N_{g.d} $ to be
the (zero-dimensional) intersection product
  \[ \ev_1^* P_1 \cdot \, \cdots \, \cdot \ev_n^* P_n \in \ZZ \]
on $ \bar M_{g,n}(r,d) $. Taking this as a definition has the advantage that we
do not have to care about whether (or for which collections of points $ P_i $)
the number of curves through the $ P_i $ actually \emph {is} finite --- we get
a well-defined number in any case. These numbers are called the Gromov-Witten
invariants.

In sections \ref {wdvv} and \ref {caporaso} we will explain how these numbers
can actually be computed --- both in Gromov-Witten theory and in tropical
geometry. For the moment let us just explain how our problem can be set up in
the tropical world. In the same way as at the end of section \ref {bezout} the
idea is of course simply to map the whole situation to the real plane by the
logarithm resp.\ valuation map. A plane curve over $ \CC $ resp.\ the field of
Puiseux series $K$ through some points $ P_1,\dots,P_n $ then simply maps to a
plane tropical curve in $ \RR^2 $ through the $n$ image points. It should
therefore also be possible to compute the numbers $ N_{g.d} $ by counting plane
\emph {tropical} curves of the given genus and degree through $n$ given points
in the real plane. The following picture shows the simplest example of this
statement, namely that also in the tropical world there is always exactly one
line through two given (general) points. Note that the relative position of the
points in the plane determines on which edges of the tropical line the points
lie (i.e.\ whether we are in case (a), (b), or (c)):

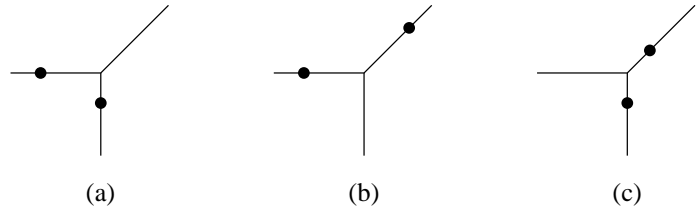
\begin {figure}[H]
  \begin {center} \input {pics/line2} \end {center}
  \caption {There is always one tropical line through two given (general)
    points}
  \label {pic-line2}
\end {figure}

This program has been carried out successfully for all $ N_{g,d} $ (and in fact
for curves in any toric surface) by Mikhalkin \cite {M}. One of the main
problems when transferring the situation to the tropical world is to determine
--- in the same way as for B\'ezout's theorem at the end of section \ref
{bezout} --- \emph {how many} complex curves through the given points map to
the same tropical curve, i.e.\ with what multiplicity the tropical curves have
to be counted. The answer to this question turns out to be surprisingly simple
and even independent of the chosen points: let $ \Gamma \subset \RR^2 $ be a
tropical curve through the given points. If the points are in general position
then all vertices of $ \Gamma $ will have valence 3. For any such vertex $V$ we
first define its multiplicity to be $ w_1 w_2 | \det (u_1,u_2) | $, where $
w_1,w_2,w_3 $ and $ u_1,u_2,u_3 $ are the weights and primitive integer vectors
along the three edges adjacent to $V$ (the balancing condition ensures that it
does not matter which two of the edges we use in the formula). The
multiplicity of $ \Gamma $ is then simply the product of the multiplicities of
all its vertices. For example, the multiplicity of the following curve is 4
(the multiplicity of both vertices is 2):

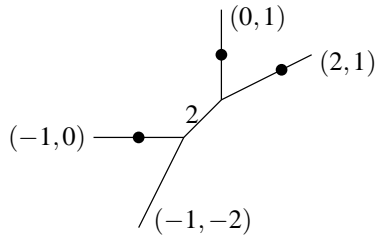
\begin {figure}[H]
  \begin {center} \input {pics/genpos1} \end {center}
  \caption {A tropical curve whose multiplicity is 4}
  \label {pic-genpos1}
\end {figure}

(Note that this is not a plane curve of some degree since the ends do not point
in the right directions for this. We can interpret figure \ref {pic-genpos1}
either as a tropical curve in a different toric surface or as a local picture
of a plane tropical curve.)

Mikhalkin's ``Correspondence Theorem'' now states that this is precisely the
correct multiplicity for our purposes, i.e.\ that the numbers $ N_{g,d} $ of
complex curves through given points $ P_1,\dots,P_n $ are the same as the
numbers of tropical curves of the same genus and degree through the images of $
P_1,\dots,P_n $ under the logarithm (resp.\ valuation) map.

One easy corollary of this statement is worth mentioning: as the numbers of
complex curves do not depend on the choice of points $ P_i $ (as long as they
are in general position) the same must be true in the tropical setting. To see
that this is in fact a non-trivial statement let us consider again the tropical
curve of figure \ref {pic-genpos1}. Because of the balancing condition this
curve is fixed in the plane by the directions of the outer edges and the
positions of the three marked points in the plane. If we now move the rightmost
point down this has the effect of shrinking the bounded edge of the curve (see
picture (a) below) until we reach a curve with a 4-valent vertex in (b) (note
that the slopes of all edges are fixed):

\begin {figure}[H]
  \begin {center} \input {pics/genpos2} \end {center}
  \caption {The weighted number of tropical curves does not change when moving
    the points}
  \label {pic-genpos2}
\end {figure}
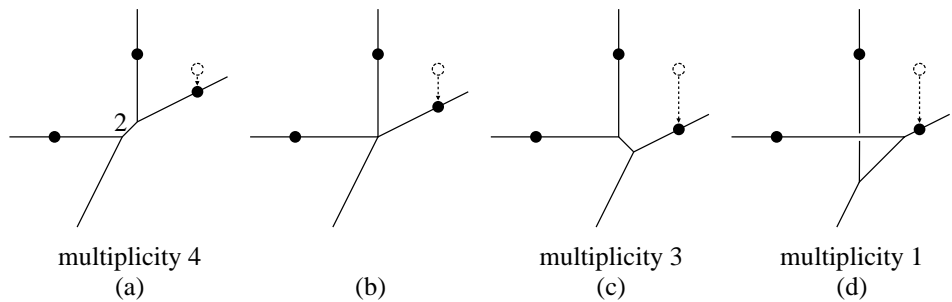

If we move the marked point down further we see that there are now \emph {two}
combinatorially different possibilities (c) and (d) what the curve might look
like. These two types have multiplicities 3 and 1 respectively, so their total
\emph {weighted} sum is the same as before, in accordance with Mikhalkin's
theorem. A purely tropical proof of the statement that the weighted number of
tropical curves through some given points does not depend on the choice of
points can be found in \cite {GM1}.

It would of course be desirable to set up tropical moduli spaces of stable maps
that are themselves tropical varieties, establish a tropical intersection
theory and use this to explain both the definition of the multiplicity of
the curves and the constancy of the numbers $ N_{g,d} $ by general principles
of this theory --- in the same way as in the complex case. So far it is not
known how to do this however.

Before we explain in the next sections how the numbers $ N_{g,d} $ can be
computed in Gromov-Witten theory and how these methods can be transferred to
the tropical setting we should mention that Mikhalkin has also found a
different way to calculate these numbers tropically by an algorithm that has no
analogue in classical geometry \cite {M}. His idea is to identify the curves in
question by their Newton subdivisions and thus translate the problem into one
of counting subdivisions of the polytope $ \Delta_d $ with certain properties.
In comparison to the algorithms that we will describe next this method has
a very high combinatorial complexity however and is therefore probably not
suited very well for actual numerical computations.

\subsection {The tropical WDVV equations} \label {wdvv}

After having discussed the tropical way to set up enumerative problems let us
now focus on how the numbers can be computed. We will first deal with \emph
{rational} curves, i.e.\ with the numbers $ N_d := N_{0,d} $ of plane curves of
genus 0 that pass through $ 3d-1 $ points. These numbers have first been
computed by Kontsevich about ten years ago. His result was that the numbers $
N_d $ are given recursively by the initial value $ N_1 = 1 $ and the equation
  \[ N_d = \sum_{\substack {d_1+d_2=d \\ d_1,d_2>0}} \left(
       d_1^2 d_2^2 \binom {3d-4}{3d_1-2} - d_1^3 d_2 \binom {3d-4}{3d_1-1}
     \right) N_{d_1} N_{d_2} \]
for $ d>1 $ (see \cite {KM} claim 5.2.1).

The main tool in deriving this formula is the so-called WDVV equations. To
explain the origin of these equations we have to study the moduli spaces of
stable maps $ \bar M_{0,n} (r,d) $ in a little more detail first. One of the
key ideas in the construction of these spaces is that it is usually better to
parametrize curves in $ \PP^r $ as maps from an abstract curve to $ \PP^r $
rather than embedded curves in $ \PP^r $ (hence the name ``stable maps''). More
precisely, the points of $ \bar M_{0,n} (r,d) $ are in bijection to tuples $
(C,x_1,\dots,x_n,f) $ where $ x_1,\dots,x_n $ are distinct smooth points on a
rational nodal curve $C$ and $ f:C \to \PP^r $ is a morphism of degree $d$
(with a stability condition). A special case of this construction is $ r=0 $
when there is no map and we simply parametrize certain ``stable'' nodal
rational curves with $n$ marked points. This gives rise to the well-known \emph
{moduli spaces of stable curves} $ \bar M_{0,n} $ that have first been
considered by Deligne and Mumford \cite {DM}. The most important example for
our purposes is the space $ \bar M_{0,4} $ of 4-pointed stable rational curves.
It is well-known that this space is isomorphic to $ \PP^1 $. The general point
of $ \bar M_{0,4} $ corresponds to a smooth rational curve with 4 distinct
marked points, whereas there are also three special points corresponding to the
curves

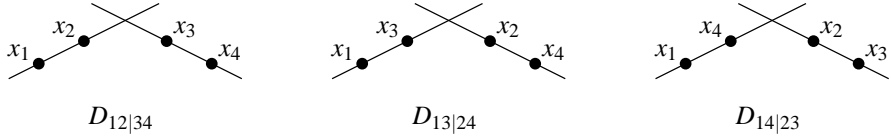
\begin {figure}[H]
  \begin {center} \input {pics/m04-points} \end {center}
  \caption {The three special points in $ \bar M_{0,4} $}
  \label {pic-m04-points}
\end {figure}

As $ \bar M_{0,4} $ is isomorphic to $ \PP^1 $ these three points obviously
define the same homology class (resp.\ the same divisor) on this moduli space.

The important point is now that there are ``forgetful maps'' $ \pi: \bar
M_{0,n} (r,d) \to \bar M_{0,4} $ for all $ n \ge 4 $ that send a stable map $
(C,x_1,\dots,x_n,f) $ to (the stabilization of) $ (C,x_1,\dots,x_4) $. Pulling
back the equality $ D_{12|34} = D_{13|24} $ of homology classes on $ \bar
M_{0,4} $ we conclude that $ \pi^* D_{12|34} $ and $ \pi^* D_{13|24} $ define
the same homology class in $ \bar M_{0,n} (r,d) $. Now $ \pi^* D_{12|34} $ (and
similarly of course for $ \pi^* D_{13|24} $) can be described explicitly
as the locus of all \emph {reducible} stable maps with two components such that
the marked points $ x_1,x_2 $ lie on one component and $ x_3,x_4 $ on the
other. Note that this space has many irreducible components since the degree
and the marked points $ x_5,\dots,x_n $ can be distributed onto the two
components in an arbitrary way.

If we now intersect the equation $ \pi^* D_{12|34} = \pi^* D_{13|24} $ of
codimension-1 cycles with suitable cycles of dimension 1 (that correspond to
the condition that the curves in question pass through given subspaces of $
\PP^r $ at the marked points) we get some equations between certain numbers of
reducible curves through given points. But these numbers of reducible curves
are just products of the corresponding numbers for their irreducible
components, i.e.\ products of certain numbers of curves of smaller degree. This
way one obtains recursion formulas that can be shown to determine all the
numbers completely (for more details see e.g.\ \cite {CK} section 7.4.2); in
the case of $ \PP^2 $ we get Kontsevich's formula stated above.

It has been shown very recently that Kontsevich's formula can also be proven in
essentially the same way in tropical geometry \cite {GM3}. For a tropical
version of the above proof it is very important that we also adapt the ``stable
map picture'' to the tropical setting, i.e.\ parametrize plane tropical curves
as maps from an ``abstract tropical curve'' to $ \PP^r $. Here by abstract
tropical curve we simply mean a connected graph $ \Gamma $ obtained by glueing
closed (not necessarily bounded) real intervals together at their boundary
points in such a way that every vertex has valence at least 3. In particular,
every bounded edge of such an abstract tropical curve has an intrinsic length.
Following an idea of Mikhalkin \cite {M3} the unbounded ends of $ \Gamma $
will be labeled and called the marked points of the curve. The most important
example for our applications is of course the ``tropical $ \bar M_{0,4} $''
whose points correspond to tree graphs with 4 unbounded ends. There are four
possible combinatorial types for this:

\begin {figure}[H]
  \begin {center} \input {pics/m04-trop} \end {center}
  \caption {The four combinatorial types of curves in the tropical $
    \bar M_{0,4} $}
  \label {pic-m04-trop}
\end {figure}
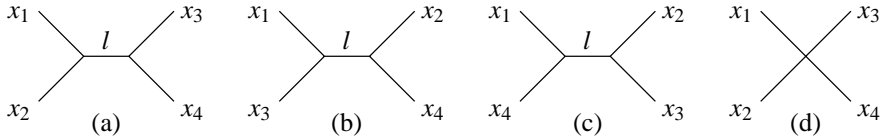

In the types (a) to (c) the bounded edge has an intrinsic length $l$; so each
of these types leads to a stratum of $ \bar M_{0,4} $ isomorphic to $ \RR_{>0}
$ parametrized by this length. The last type (d) is simply a point in $
\bar M_{0,4} $ that can be seen as the boundary point where the other three
strata meet. Hence $ \bar M_{0,4} $ is again a rational tropical curve:

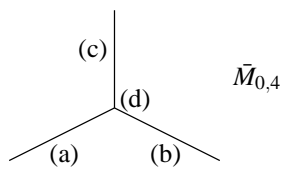
\begin {figure}[H]
  \begin {center} \input {pics/m04-trop-mod} \end {center}
  \caption {$ \bar M_{0,4} $ is again a rational tropical curve}
  \label {pic-m04-trop-mod}
\end {figure}

In analogy to the complex case plane tropical curves are now parametrized as
tuples $ (\Gamma,x_1,\dots,x_n,h) $, where $ \Gamma $ is an abstract tropical
curve, $ x_1,\dots,x_n $ are distinct unbounded ends of $ \Gamma $, and $ h:
\Gamma \to \RR^2 $ is a piecewise linear map with certain conditions (see \cite
{GM3} for details). The most important feature of this definition is that $h$
may be a constant map on some edges of $ \Gamma $, and is in fact \emph
{required} to be a constant map on the unbounded ends $ x_1,\dots,x_n $.  For
example, the following picture shows a 4-pointed plane tropical conic, i.e.\ of
the tropical analogue of $ \bar M_{0,4} (2,2) $:

\begin {figure}[H]
  \begin {center} \input {pics/trop-conic-irr} \end {center}
  \caption {A tropical stable map in $ \bar M_{0,4}(2,2) $}
  \label {pic-trop-conic-irr}
\end {figure}
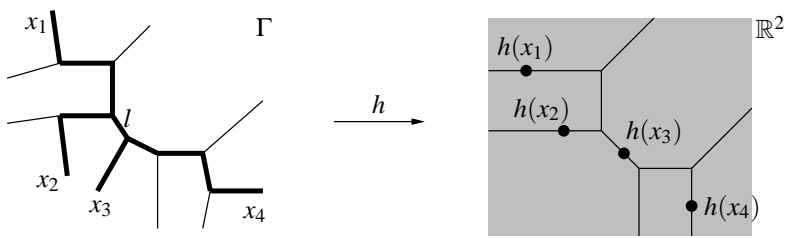

Note that the balancing condition around every (3-valent) vertex adjacent to a
marked point (resp.\ edge) $ x_i $ ensures that the other two edges around this
vertex form a straight line in $ \RR^2 $.

It is easy to see from this picture already that the tropical moduli spaces $
\bar M_{0,n}(r,d) $ with $ n \ge 4 $ admit forgetful maps to $ \bar M_{0,4} $:
given an $n$-marked plane tropical curve $ (\Gamma,x_1,\dots,x_n,h) $ we simply
forget the map $h$, take the minimal connected subgraph of $ \Gamma $ that
contains $ x_1,\dots,x_4 $, and ``straighten'' this graph to obtain an element
of $ \bar M_{0,4} $. In the picture above we simply obtain the ``straightened
version'' of the subgraph drawn in bold, i.e.\ the element of $ \bar M_{0,4} $
of type (a) in figure \ref {pic-m04-trop} with length parameter $l$ as
indicated in the picture.

To obtain the WDVV equations we now simply consider the inverse image under
this forgetful map of a point of $ \bar M_{0,4} $ of type (a) resp.\ (b) in
figure \ref {pic-m04-trop} with a very large length parameter $l$. It can
be shown that such very large lengths can occur only if there is a bounded edge
(of a very large length) in $ \Gamma $ on which $h$ is constant:

\begin {figure}[H]
  \begin {center} \input {pics/trop-conic-red} \end {center}
  \caption {A very large length in $ \bar M_{0,4} $ leads to a tropical stable
    map with a contracted bounded edge}
  \label {pic-trop-conic-red}
\end {figure}
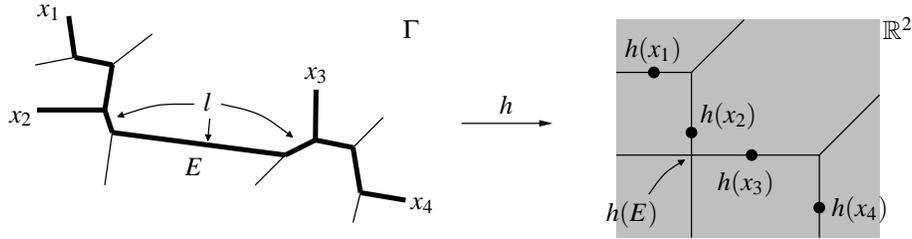

Again the balancing condition around the contracted bounded edge requires that
the image of the tropical curve in $ \RR^2 $ is locally a union of two straight
lines. We can therefore consider such curves as being reducible and made up of
two tropical curves of smaller degrees (in figure \ref {pic-trop-conic-red} we
have a reducible tropical conic that is a union of two tropical lines). The
picture is now exactly the same as in the classical case, and in fact the rest
of the proof of Kontsevich's formula works in the same way as in Gromov-Witten
theory. It is expected that essentially the same proof can be used to reprove
the WDVV equations for rational curves in higher-dimensional spaces as well.

This application of tropical geometry shows very well that it should be
possible to carry many concepts from classical complex geometry over to the
tropical world: moduli spaces of curves and stable maps, morphisms, divisors
and divisor classes, intersection multiplicities, and so on. In \cite {GM3}
these concepts were introduced only in the specific cases needed for
Kontsevich's formula.

\subsection {The tropical Caporaso-Harris formula} \label {caporaso}

After having discussed rational curves let us now turn to the general numbers $
N_{g,d} $ for arbitrary genus $g$. These numbers have first been computed by
Caporaso and Harris \cite {CH}. The idea in their proof is is to define new
invariants that count plane curves of given degree and genus having specified
local contact orders to a fixed line $L$ and passing in addition through the
appropriate number of general points. By specializing one point after the other
to lie on $L$ one can then derive recursive relations among these new
invariants that finally suffice to compute all the numbers $ N_{g,d} $.

Instead of explaining the general formula let us look at an example of what
happens in this specialization process, referring to \cite {CH} for details. We
consider plane rational cubics having a point of contact order 3 to $L$ at a
fixed point $ P_1 \in L $ and passing in addition through 5 general points $
P_2,\dots,P_6 \in \PP^2 $ as in the following picture on the left:

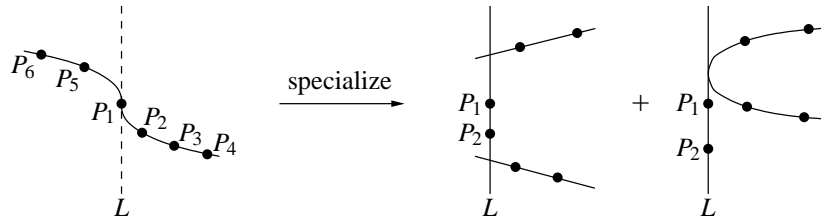
\begin {figure}[H]
  \begin {center} \input {pics/ch-complex} \end {center}
  \caption {Computing the number of cubics with a point of contact order 3 to a
    line}
  \label {pic-ch-complex}
\end {figure}

To compute the number of such curves we specialize $ P_2 $ to lie on $L$. As
the cubics intersect $L$ already with multiplicity 3 at $ P_1 $ they cannot
pass through another point on $L$ unless they become reducible and have $L$ as
a component. Hence there are two possibilities after the specialization (see
figure \ref {pic-ch-complex}): the cubics can degenerate into a union of three
lines $ L \cup L_1 \cup L_2 $ where $ L_1 $ and $ L_2 $ each pass through two
of the points $ P_3,\dots,P_6 $, or they can degenerate into $ L \cup C $,
where $C$ is a conic tangent to $L$ and passing through $ P_3,\dots,P_6 $.
The initial number of rational cubics with a point of contact order 3 to $L$ at
a fixed point and passing through 5 more general points is therefore a sum of
two numbers (counted with suitable multiplicities) related to only lines and
conics. This is the general idea of Caporaso and Harris how specialization
finally reduces the degree of the curves and allows a recursive solution to
compute the numbers $ N_{g,d} $ as well as all the newly introduced numbers of
curves with multiplicity conditions.

In fact, the same constructions can again be made in tropical geometry \cite
{GM2}. Intuitively, if we pick our line $ L \subset \CC^2 $ to be the line with
$ z_1 $-coordinate 0 then its image under the logarithm map is ``the vertical
line with $ x_1 $-coordinate $ -\infty $. Hence the process of moving $ P_2 $
to the line $L$ in complex geometry now simply corresponds to moving $ P_2 $ to
the very far left in tropical geometry. Moreover, curves with higher contact
orders to $L$ in complex geometry just correspond to tropical curves with
unbounded ends of higher weight to the left. So the tropical analogue of the
specialization process of figure \ref {pic-ch-complex} is

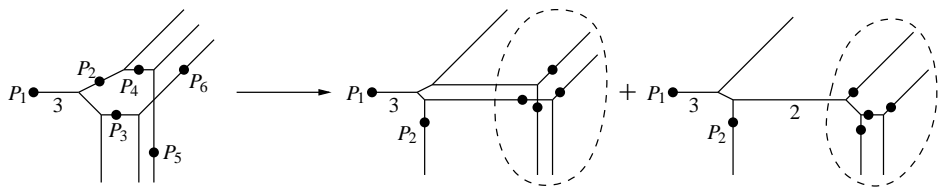
\begin {figure}[H]
  \begin {center} \input {pics/ch-trop} \end {center}
  \caption {The degeneration of figure \ref {pic-ch-complex} in the tropical
    setting}
  \label {pic-ch-trop}
\end {figure}

where $ P_1 $ is to be considered to lie infinitely far to the left. Note that
the curves after moving $ P_2 $ to the left are not reducible, but they still
``split'' into two parts: a left part (through $ P_2 $) and a right part
(through the remaining points, circled in the picture above). We get the same
``degenerations'' as in the complex case: one where the right part consists of
two lines through two of the points $ P_3,\dots,P_6 $ each, and one where it
consists of a conic ``tangent to a line'' (i.e.\ with an unbounded edge of
multiplicity 2 to the left).

Using this idea it has been shown in \cite {GM2} that the Caporaso-Harris
formula can also be proven in tropical geometry --- in fact with a much simpler
proof than in complex geometry since we are only dealing with combinatorial
objects and do not have to construct and study complicated moduli spaces in
complex geometry.

\subsection {Real enumerative geometry and Welschinger invariants}
  \label {welschinger}

Of course, the same questions as in the previous sections can be asked for real
instead of for complex curves: given $ d \ge 1 $, $ g \ge 0 $, and $ 3d+g-1 $
points in general position in the real projective plane $ \PP^2_\RR $, how many
curves of degree $d$ and genus $g$ are there in $ \PP^2_\RR $ that intersect
all the given points? Note that every such real curve has a complexification
by just considering its equation as an equation in the complex rather than the
real plane, and by its degree and genus we simply mean the degree resp.\ genus
of its complexification.

As usual in algebraic geometry this real case is much more difficult to handle
than the complex case. The first problem is already that the answer to this
question will in general depend on the position of the points so that there are
no well-defined numbers $ N_{g,d} $ as in the complex case. Instead one could
ask questions of the following type:
\begin {itemize}
\item is there at least one real curve of degree $d$ and genus $g$ through any
  choice of $ 3d+g-1 $ given points? Or even better: can we compute a lower
  bound for the number of real curves through any choice of points?
\item Is there a way to assign multiplicities to the real curves through the
  given points so that the \emph {weighted} sum of these curves is independent
  of the choice of points?
\end {itemize}
A few years ago Welschinger has found a solution to the second question for the
case of curves of genus 0 \cite {W}. To state his result note first that a
general \emph {complex} plane rational curve of degree $d$ has exactly $ \frac
12(d-1)(d-2) $ nodes, i.e.\ points where two smooth branches of the curve
intersect transversally. If $C$ is now a real rational plane curve then each
of the nodes of its complexification is of one of the following three types:
\begin {enumerate}
\item nodes that are in the real plane $ \PP^2_\RR $ and where the local
  equation of the curve is of the form $ x^2-y^2=0 $, i.e.\ $ (x+y)(x-y) = 0 $
  for suitable local analytic coordinates. In a local real picture $C$ is
  simply the union of two smooth curves intersecting transversely.
\item nodes that are in the real plane $ \PP^2_\RR $ and where the local
  equation of the curve is of the form $ x^2+y^2=0 $ for suitable local
  analytic coordinates. In a local real picture such a node leads to an
  isolated point corresponding to the values $ x=y=0 $.
\item nodes that are not in the real plane $ \PP^2_\RR $. These nodes obviously
  come in pairs as the complex conjugate of such a node is again a node of the
  complexification. They are not visible in the real curve $C$.
\end {enumerate}
Welschinger's main theorem is now the following: if we assign to each real
plane curve $C$ the multiplicity $ (-1)^m $ where $m$ is the number of nodes of
type (b) of its complexification then the corresponding weighted sum of all
curves through the given points is independent of the choice of points.
It is called the \emph {Welschinger invariant} $ W_d $. In particular, $ W_d $
gives a lower bound for the actual number of real plane curves through any set
of $ 3d-1 $ general given points, whereas the complex number $ N_{0,d} $ is of
course an upper bound.

Unfortunately, except for a few special cases Welschinger was not able to
actually compute the numbers $ W_d $. So the question whether there always
exists a real rational plane curve of degree $d$ through any set of $ 3d-1 $
general points remained open.

Some time ago however a tropical way has been found to compute the Welschinger
invariants $ W_d $ and in fact also the actual (non-invariant) numbers of real
curves through some configurations of points in the plane \cite {IKS,M}. In the
same way as in Mikhalkin's Correspondence Theorem the strategy is to identify
complex resp.\ real algebraic curves with tropical curves, and then to count
such tropical curves with the proper multiplicities. More precisely, to obtain
the Welschinger invariant $ W_d $ one has to count rational plane tropical
curves of degree $d$ through $ 3d-1 $ points in the same way as for the
computation of $ N_{0,d} $ --- but one does not count them with the
``complex multiplicity'' as described in section \ref {gromov-witten}, but
rather with the multiplicity
\begin {align*}
  0 & \quad \mbox {if the ``complex multiplicity'' is even}, \\
  1 & \quad \mbox {if the ``complex multiplicity'' is congruent to 1 modulo
    4}, \\
  -1 & \quad \mbox {if the ``complex multiplicity'' is congruent to 3 modulo
    4}.
\end {align*}
One can then try to do this count by enumerating the corresponding Newton
subdivisions of the polytope $ \Delta_d $ as in section \ref {graphs}. Although
the algorithm is combinatorially very complicated (and cannot be explained here
in detail) it can be used to prove that the numbers $ W_d $ are all positive
and hence that there is always at least one real rational plane curve of degree
$d$ through any set of $ 3d-1 $ points in general position. In fact, a more
careful study of the algorithm even allows one to prove that the lower and
upper bounds $ W_d $ resp.\ $ N_{0,d} $ for the numbers of these curves grow
approximately with the same speed as $d$ increases \cite {IKS2}.

Very recently Itenberg has shown that both the proof of \cite {GM1} that the
invariants are independent of the marked points and the proof of \cite {GM2} of
the tropical Caporaso-Harris formula (see section \ref {caporaso}) can
be adapted to the Welschinger case. In particular, this yields a ``real
Caporaso-Harris formula'' that gives a fast method to compute the Welschinger
invariants. It is not known yet whether there exists an analogue of the WDVV
equations (see section \ref {wdvv}) in the real case.

%% file: pics/line2.tex
\begin{picture}(0,0)%
\includegraphics{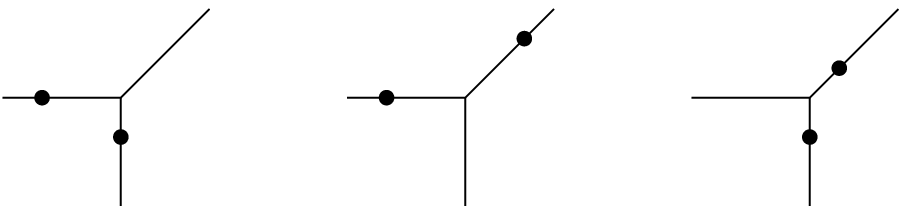}%
\end{picture}%
\setlength{\unitlength}{4144sp}%
\begingroup\makeatletter\ifx\SetFigFont\undefined%
\gdef\SetFigFont#1#2#3#4#5{%
  \reset@font\fontsize{#1}{#2pt}%
  \fontfamily{#3}\fontseries{#4}\fontshape{#5}%
  \selectfont}%
\fi\endgroup%
\begin{picture}(4119,1237)(484,-836)
\put(1036,-781){\makebox(0,0)[b]{\smash{{\SetFigFont{10}{12.0}{\familydefault}{\mddefault}{\updefault}{\color[rgb]{0,0,0}(a)}%
}}}}
\put(2611,-781){\makebox(0,0)[b]{\smash{{\SetFigFont{10}{12.0}{\familydefault}{\mddefault}{\updefault}{\color[rgb]{0,0,0}(b)}%
}}}}
\put(4186,-781){\makebox(0,0)[b]{\smash{{\SetFigFont{10}{12.0}{\familydefault}{\mddefault}{\updefault}{\color[rgb]{0,0,0}(c)}%
}}}}
\end{picture}%

%% file: pics/genpos1.tex
\begin{picture}(0,0)%
\includegraphics{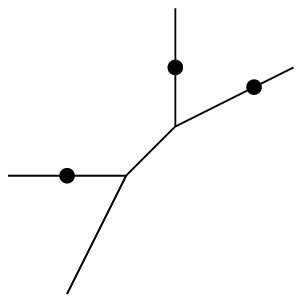}%
\end{picture}%
\setlength{\unitlength}{4144sp}%
\begingroup\makeatletter\ifx\SetFigFont\undefined%
\gdef\SetFigFont#1#2#3#4#5{%
  \reset@font\fontsize{#1}{#2pt}%
  \fontfamily{#3}\fontseries{#4}\fontshape{#5}%
  \selectfont}%
\fi\endgroup%
\begin{picture}(1425,1410)(391,-661)
\put(406,-106){\makebox(0,0)[rb]{\smash{{\SetFigFont{10}{12.0}{\familydefault}{\mddefault}{\updefault}{\color[rgb]{0,0,0}$(-1,0)$}%
}}}}
\put(1081, 29){\makebox(0,0)[rb]{\smash{{\SetFigFont{10}{12.0}{\familydefault}{\mddefault}{\updefault}{\color[rgb]{0,0,0}$2$}%
}}}}
\put(1261,614){\makebox(0,0)[lb]{\smash{{\SetFigFont{10}{12.0}{\familydefault}{\mddefault}{\updefault}{\color[rgb]{0,0,0}$(0,1)$}%
}}}}
\put(1801,344){\makebox(0,0)[lb]{\smash{{\SetFigFont{10}{12.0}{\familydefault}{\mddefault}{\updefault}{\color[rgb]{0,0,0}$(2,1)$}%
}}}}
\put(811,-601){\makebox(0,0)[lb]{\smash{{\SetFigFont{10}{12.0}{\familydefault}{\mddefault}{\updefault}{\color[rgb]{0,0,0}$(-1,-2)$}%
}}}}
\end{picture}%

%% file: pics/genpos2.tex
\begin{picture}(0,0)%
\includegraphics{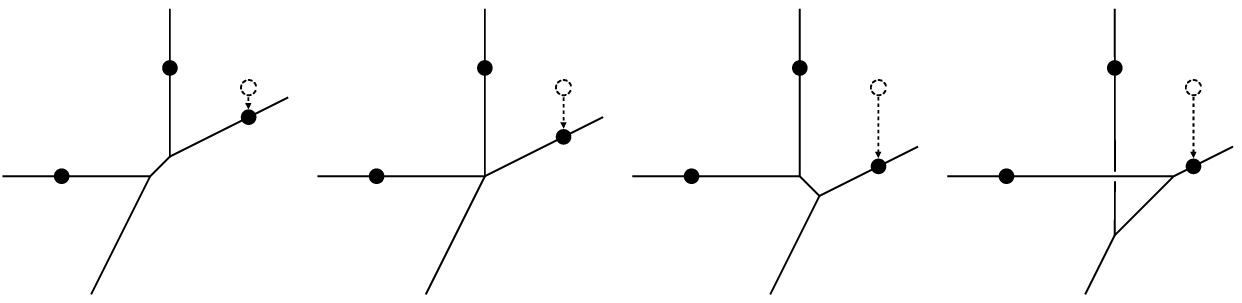}%
\end{picture}%
\setlength{\unitlength}{4144sp}%
\begingroup\makeatletter\ifx\SetFigFont\undefined%
\gdef\SetFigFont#1#2#3#4#5{%
  \reset@font\fontsize{#1}{#2pt}%
  \fontfamily{#3}\fontseries{#4}\fontshape{#5}%
  \selectfont}%
\fi\endgroup%
\begin{picture}(5649,1777)(439,-1061)
\put(1156,-26){\makebox(0,0)[rb]{\smash{{\SetFigFont{10}{12.0}{\familydefault}{\mddefault}{\updefault}{\color[rgb]{0,0,0}$2$}%
}}}}
\put(1171,-1006){\makebox(0,0)[b]{\smash{{\SetFigFont{10}{12.0}{\familydefault}{\mddefault}{\updefault}{\color[rgb]{0,0,0}(a)}%
}}}}
\put(2611,-1006){\makebox(0,0)[b]{\smash{{\SetFigFont{10}{12.0}{\familydefault}{\mddefault}{\updefault}{\color[rgb]{0,0,0}(b)}%
}}}}
\put(1171,-826){\makebox(0,0)[b]{\smash{{\SetFigFont{10}{12.0}{\familydefault}{\mddefault}{\updefault}{\color[rgb]{0,0,0}multiplicity 4}%
}}}}
\put(4051,-1006){\makebox(0,0)[b]{\smash{{\SetFigFont{10}{12.0}{\familydefault}{\mddefault}{\updefault}{\color[rgb]{0,0,0}(c)}%
}}}}
\put(4051,-826){\makebox(0,0)[b]{\smash{{\SetFigFont{10}{12.0}{\familydefault}{\mddefault}{\updefault}{\color[rgb]{0,0,0}multiplicity 3}%
}}}}
\put(5491,-826){\makebox(0,0)[b]{\smash{{\SetFigFont{10}{12.0}{\familydefault}{\mddefault}{\updefault}{\color[rgb]{0,0,0}multiplicity 1}%
}}}}
\put(5491,-1006){\makebox(0,0)[b]{\smash{{\SetFigFont{10}{12.0}{\familydefault}{\mddefault}{\updefault}{\color[rgb]{0,0,0}(d)}%
}}}}
\end{picture}%

%% file: pics/m04-points.tex
\begin{picture}(0,0)%
\includegraphics{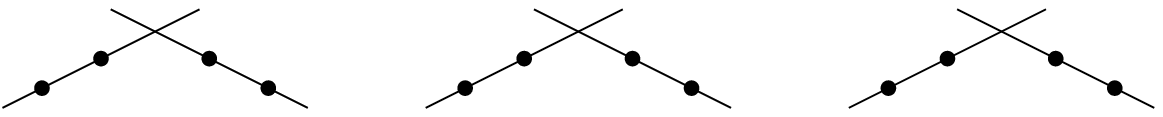}%
\end{picture}%
\setlength{\unitlength}{4144sp}%
\begingroup\makeatletter\ifx\SetFigFont\undefined%
\gdef\SetFigFont#1#2#3#4#5{%
  \reset@font\fontsize{#1}{#2pt}%
  \fontfamily{#3}\fontseries{#4}\fontshape{#5}%
  \selectfont}%
\fi\endgroup%
\begin{picture}(5289,796)(439,-440)
\put(1441,164){\makebox(0,0)[lb]{\smash{{\SetFigFont{10}{12.0}{\familydefault}{\mddefault}{\updefault}{\color[rgb]{0,0,0}$x_3$}%
}}}}
\put(1711, 29){\makebox(0,0)[lb]{\smash{{\SetFigFont{10}{12.0}{\familydefault}{\mddefault}{\updefault}{\color[rgb]{0,0,0}$x_4$}%
}}}}
\put(856,164){\makebox(0,0)[rb]{\smash{{\SetFigFont{10}{12.0}{\familydefault}{\mddefault}{\updefault}{\color[rgb]{0,0,0}$x_2$}%
}}}}
\put(586, 29){\makebox(0,0)[rb]{\smash{{\SetFigFont{10}{12.0}{\familydefault}{\mddefault}{\updefault}{\color[rgb]{0,0,0}$x_1$}%
}}}}
\put(3646, 29){\makebox(0,0)[lb]{\smash{{\SetFigFont{10}{12.0}{\familydefault}{\mddefault}{\updefault}{\color[rgb]{0,0,0}$x_4$}%
}}}}
\put(2521, 29){\makebox(0,0)[rb]{\smash{{\SetFigFont{10}{12.0}{\familydefault}{\mddefault}{\updefault}{\color[rgb]{0,0,0}$x_1$}%
}}}}
\put(2791,164){\makebox(0,0)[rb]{\smash{{\SetFigFont{10}{12.0}{\familydefault}{\mddefault}{\updefault}{\color[rgb]{0,0,0}$x_3$}%
}}}}
\put(3376,164){\makebox(0,0)[lb]{\smash{{\SetFigFont{10}{12.0}{\familydefault}{\mddefault}{\updefault}{\color[rgb]{0,0,0}$x_2$}%
}}}}
\put(4456, 29){\makebox(0,0)[rb]{\smash{{\SetFigFont{10}{12.0}{\familydefault}{\mddefault}{\updefault}{\color[rgb]{0,0,0}$x_1$}%
}}}}
\put(5311,164){\makebox(0,0)[lb]{\smash{{\SetFigFont{10}{12.0}{\familydefault}{\mddefault}{\updefault}{\color[rgb]{0,0,0}$x_2$}%
}}}}
\put(1126,-376){\makebox(0,0)[b]{\smash{{\SetFigFont{10}{12.0}{\familydefault}{\mddefault}{\updefault}{\color[rgb]{0,0,0}$ D_{12|34} $}%
}}}}
\put(3061,-376){\makebox(0,0)[b]{\smash{{\SetFigFont{10}{12.0}{\familydefault}{\mddefault}{\updefault}{\color[rgb]{0,0,0}$ D_{13|24} $}%
}}}}
\put(4726,164){\makebox(0,0)[rb]{\smash{{\SetFigFont{10}{12.0}{\familydefault}{\mddefault}{\updefault}{\color[rgb]{0,0,0}$x_4$}%
}}}}
\put(5581, 29){\makebox(0,0)[lb]{\smash{{\SetFigFont{10}{12.0}{\familydefault}{\mddefault}{\updefault}{\color[rgb]{0,0,0}$x_3$}%
}}}}
\put(4996,-376){\makebox(0,0)[b]{\smash{{\SetFigFont{10}{12.0}{\familydefault}{\mddefault}{\updefault}{\color[rgb]{0,0,0}$ D_{14|23} $}%
}}}}
\end{picture}%

%% file: pics/m04-trop.tex
\begin{picture}(0,0)%
\includegraphics{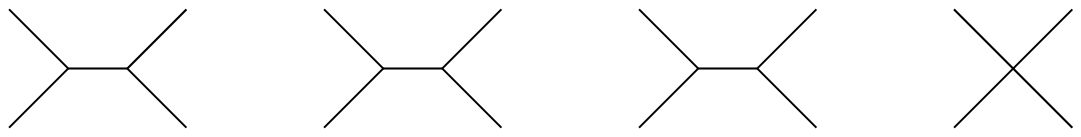}%
\end{picture}%
\setlength{\unitlength}{4144sp}%
\begingroup\makeatletter\ifx\SetFigFont\undefined%
\gdef\SetFigFont#1#2#3#4#5{%
  \reset@font\fontsize{#1}{#2pt}%
  \fontfamily{#3}\fontseries{#4}\fontshape{#5}%
  \selectfont}%
\fi\endgroup%
\begin{picture}(4980,877)(346,-296)
\put(1261,434){\makebox(0,0)[lb]{\smash{{\SetFigFont{10}{12.0}{\familydefault}{\mddefault}{\updefault}{\color[rgb]{0,0,0}$x_3$}%
}}}}
\put(361,434){\makebox(0,0)[rb]{\smash{{\SetFigFont{10}{12.0}{\familydefault}{\mddefault}{\updefault}{\color[rgb]{0,0,0}$x_1$}%
}}}}
\put(361,-151){\makebox(0,0)[rb]{\smash{{\SetFigFont{10}{12.0}{\familydefault}{\mddefault}{\updefault}{\color[rgb]{0,0,0}$x_2$}%
}}}}
\put(1801,434){\makebox(0,0)[rb]{\smash{{\SetFigFont{10}{12.0}{\familydefault}{\mddefault}{\updefault}{\color[rgb]{0,0,0}$x_1$}%
}}}}
\put(4681,434){\makebox(0,0)[rb]{\smash{{\SetFigFont{10}{12.0}{\familydefault}{\mddefault}{\updefault}{\color[rgb]{0,0,0}$x_1$}%
}}}}
\put(4681,-151){\makebox(0,0)[rb]{\smash{{\SetFigFont{10}{12.0}{\familydefault}{\mddefault}{\updefault}{\color[rgb]{0,0,0}$x_2$}%
}}}}
\put(3241,434){\makebox(0,0)[rb]{\smash{{\SetFigFont{10}{12.0}{\familydefault}{\mddefault}{\updefault}{\color[rgb]{0,0,0}$x_1$}%
}}}}
\put(1261,-151){\makebox(0,0)[lb]{\smash{{\SetFigFont{10}{12.0}{\familydefault}{\mddefault}{\updefault}{\color[rgb]{0,0,0}$x_4$}%
}}}}
\put(2701,-151){\makebox(0,0)[lb]{\smash{{\SetFigFont{10}{12.0}{\familydefault}{\mddefault}{\updefault}{\color[rgb]{0,0,0}$x_4$}%
}}}}
\put(5311,434){\makebox(0,0)[lb]{\smash{{\SetFigFont{10}{12.0}{\familydefault}{\mddefault}{\updefault}{\color[rgb]{0,0,0}$x_3$}%
}}}}
\put(5311,-151){\makebox(0,0)[lb]{\smash{{\SetFigFont{10}{12.0}{\familydefault}{\mddefault}{\updefault}{\color[rgb]{0,0,0}$x_4$}%
}}}}
\put(1801,-151){\makebox(0,0)[rb]{\smash{{\SetFigFont{10}{12.0}{\familydefault}{\mddefault}{\updefault}{\color[rgb]{0,0,0}$x_3$}%
}}}}
\put(2701,434){\makebox(0,0)[lb]{\smash{{\SetFigFont{10}{12.0}{\familydefault}{\mddefault}{\updefault}{\color[rgb]{0,0,0}$x_2$}%
}}}}
\put(3241,-151){\makebox(0,0)[rb]{\smash{{\SetFigFont{10}{12.0}{\familydefault}{\mddefault}{\updefault}{\color[rgb]{0,0,0}$x_4$}%
}}}}
\put(4141,434){\makebox(0,0)[lb]{\smash{{\SetFigFont{10}{12.0}{\familydefault}{\mddefault}{\updefault}{\color[rgb]{0,0,0}$x_2$}%
}}}}
\put(4141,-151){\makebox(0,0)[lb]{\smash{{\SetFigFont{10}{12.0}{\familydefault}{\mddefault}{\updefault}{\color[rgb]{0,0,0}$x_3$}%
}}}}
\put(811,254){\makebox(0,0)[b]{\smash{{\SetFigFont{10}{12.0}{\familydefault}{\mddefault}{\updefault}{\color[rgb]{0,0,0}$l$}%
}}}}
\put(2251,254){\makebox(0,0)[b]{\smash{{\SetFigFont{10}{12.0}{\familydefault}{\mddefault}{\updefault}{\color[rgb]{0,0,0}$l$}%
}}}}
\put(3691,254){\makebox(0,0)[b]{\smash{{\SetFigFont{10}{12.0}{\familydefault}{\mddefault}{\updefault}{\color[rgb]{0,0,0}$l$}%
}}}}
\put(811,-241){\makebox(0,0)[b]{\smash{{\SetFigFont{10}{12.0}{\familydefault}{\mddefault}{\updefault}{\color[rgb]{0,0,0}(a)}%
}}}}
\put(2251,-241){\makebox(0,0)[b]{\smash{{\SetFigFont{10}{12.0}{\familydefault}{\mddefault}{\updefault}{\color[rgb]{0,0,0}(b)}%
}}}}
\put(3691,-241){\makebox(0,0)[b]{\smash{{\SetFigFont{10}{12.0}{\familydefault}{\mddefault}{\updefault}{\color[rgb]{0,0,0}(c)}%
}}}}
\put(4996,-241){\makebox(0,0)[b]{\smash{{\SetFigFont{10}{12.0}{\familydefault}{\mddefault}{\updefault}{\color[rgb]{0,0,0}(d)}%
}}}}
\end{picture}%

%% file: pics/m04-trop-mod.tex
\begin{picture}(0,0)%
\includegraphics{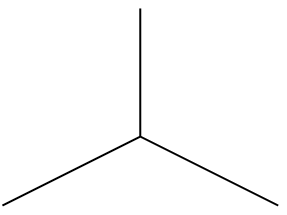}%
\end{picture}%
\setlength{\unitlength}{4144sp}%
\begingroup\makeatletter\ifx\SetFigFont\undefined%
\gdef\SetFigFont#1#2#3#4#5{%
  \reset@font\fontsize{#1}{#2pt}%
  \fontfamily{#3}\fontseries{#4}\fontshape{#5}%
  \selectfont}%
\fi\endgroup%
\begin{picture}(1512,978)(934,-1207)
\put(1276,-1147){\makebox(0,0)[b]{\smash{{\SetFigFont{10}{12.0}{\familydefault}{\mddefault}{\updefault}{\color[rgb]{0,0,0}(a)}%
}}}}
\put(1534,-523){\makebox(0,0)[rb]{\smash{{\SetFigFont{10}{12.0}{\familydefault}{\mddefault}{\updefault}{\color[rgb]{0,0,0}(c)}%
}}}}
\put(1606,-821){\makebox(0,0)[lb]{\smash{{\SetFigFont{10}{12.0}{\familydefault}{\mddefault}{\updefault}{\color[rgb]{0,0,0}(d)}%
}}}}
\put(1882,-1152){\makebox(0,0)[b]{\smash{{\SetFigFont{10}{12.0}{\familydefault}{\mddefault}{\updefault}{\color[rgb]{0,0,0}(b)}%
}}}}
\put(2431,-691){\makebox(0,0)[b]{\smash{{\SetFigFont{10}{12.0}{\familydefault}{\mddefault}{\updefault}{\color[rgb]{0,0,0}$ \bar M_{0,4} $}%
}}}}
\end{picture}%

%% file: pics/trop-conic-irr.tex
\begin{picture}(0,0)%
\includegraphics{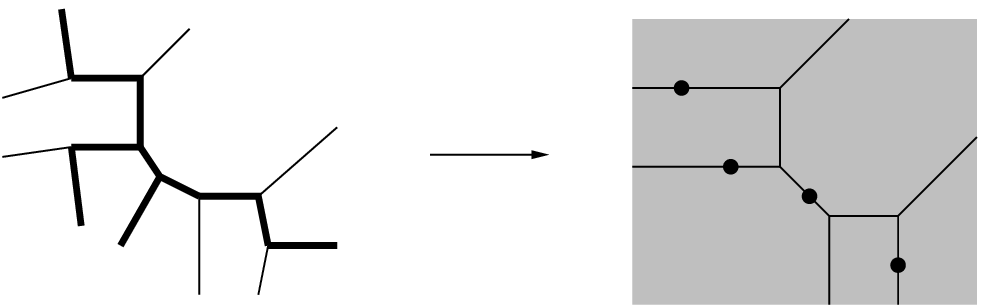}%
\end{picture}%
\setlength{\unitlength}{4144sp}%
\begingroup\makeatletter\ifx\SetFigFont\undefined%
\gdef\SetFigFont#1#2#3#4#5{%
  \reset@font\fontsize{#1}{#2pt}%
  \fontfamily{#3}\fontseries{#4}\fontshape{#5}%
  \selectfont}%
\fi\endgroup%
\begin{picture}(4502,1396)(1429,-1243)
\put(1690,  6){\makebox(0,0)[rb]{\smash{{\SetFigFont{10}{12.0}{\familydefault}{\mddefault}{\updefault}{\color[rgb]{0,0,0}$x_1$}%
}}}}
\put(1756,-939){\makebox(0,0)[rb]{\smash{{\SetFigFont{10}{12.0}{\familydefault}{\mddefault}{\updefault}{\color[rgb]{0,0,0}$x_2$}%
}}}}
\put(2002,-1086){\makebox(0,0)[b]{\smash{{\SetFigFont{10}{12.0}{\familydefault}{\mddefault}{\updefault}{\color[rgb]{0,0,0}$x_3$}%
}}}}
\put(2856,-1126){\makebox(0,0)[lb]{\smash{{\SetFigFont{10}{12.0}{\familydefault}{\mddefault}{\updefault}{\color[rgb]{0,0,0}$x_4$}%
}}}}
\put(3663,-492){\makebox(0,0)[b]{\smash{{\SetFigFont{10}{12.0}{\familydefault}{\mddefault}{\updefault}{\color[rgb]{0,0,0}$h$}%
}}}}
\put(5916,-46){\makebox(0,0)[lb]{\smash{{\SetFigFont{10}{12.0}{\familydefault}{\mddefault}{\updefault}{\color[rgb]{0,0,0}$ \RR^2 $}%
}}}}
\put(2976,-32){\makebox(0,0)[b]{\smash{{\SetFigFont{10}{12.0}{\familydefault}{\mddefault}{\updefault}{\color[rgb]{0,0,0}$ \Gamma $}%
}}}}
\put(4543,-150){\makebox(0,0)[b]{\smash{{\SetFigFont{10}{12.0}{\familydefault}{\mddefault}{\updefault}{\color[rgb]{0,0,0}$h(x_1)$}%
}}}}
\put(4810,-517){\makebox(0,0)[rb]{\smash{{\SetFigFont{10}{12.0}{\familydefault}{\mddefault}{\updefault}{\color[rgb]{0,0,0}$h(x_2)$}%
}}}}
\put(5143,-650){\makebox(0,0)[lb]{\smash{{\SetFigFont{10}{12.0}{\familydefault}{\mddefault}{\updefault}{\color[rgb]{0,0,0}$h(x_3)$}%
}}}}
\put(5610,-1104){\makebox(0,0)[lb]{\smash{{\SetFigFont{10}{12.0}{\familydefault}{\mddefault}{\updefault}{\color[rgb]{0,0,0}$h(x_4)$}%
}}}}
\put(2135,-582){\makebox(0,0)[lb]{\smash{{\SetFigFont{10}{12.0}{\familydefault}{\mddefault}{\updefault}{\color[rgb]{0,0,0}$l$}%
}}}}
\end{picture}%

%% file: pics/trop-conic-red.tex
\begin{picture}(0,0)%
\includegraphics{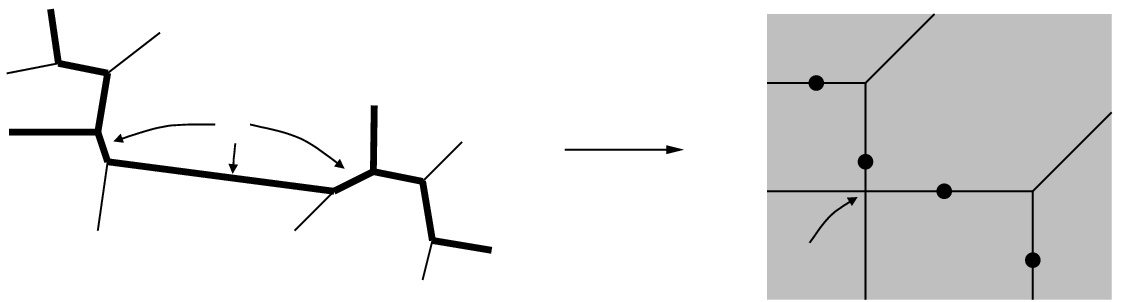}%
\end{picture}%
\setlength{\unitlength}{4144sp}%
\begingroup\makeatletter\ifx\SetFigFont\undefined%
\gdef\SetFigFont#1#2#3#4#5{%
  \reset@font\fontsize{#1}{#2pt}%
  \fontfamily{#3}\fontseries{#4}\fontshape{#5}%
  \selectfont}%
\fi\endgroup%
\begin{picture}(5117,1481)(814,-1243)
\put(3663,-492){\makebox(0,0)[b]{\smash{{\SetFigFont{10}{12.0}{\familydefault}{\mddefault}{\updefault}{\color[rgb]{0,0,0}$h$}%
}}}}
\put(5916,-46){\makebox(0,0)[lb]{\smash{{\SetFigFont{10}{12.0}{\familydefault}{\mddefault}{\updefault}{\color[rgb]{0,0,0}$ \RR^2 $}%
}}}}
\put(3092,-36){\makebox(0,0)[b]{\smash{{\SetFigFont{10}{12.0}{\familydefault}{\mddefault}{\updefault}{\color[rgb]{0,0,0}$ \Gamma $}%
}}}}
\put(999, 91){\makebox(0,0)[rb]{\smash{{\SetFigFont{10}{12.0}{\familydefault}{\mddefault}{\updefault}{\color[rgb]{0,0,0}$x_1$}%
}}}}
\put(829,-540){\makebox(0,0)[rb]{\smash{{\SetFigFont{10}{12.0}{\familydefault}{\mddefault}{\updefault}{\color[rgb]{0,0,0}$x_2$}%
}}}}
\put(4543,-151){\makebox(0,0)[b]{\smash{{\SetFigFont{10}{12.0}{\familydefault}{\mddefault}{\updefault}{\color[rgb]{0,0,0}$h(x_1)$}%
}}}}
\put(4824,-544){\makebox(0,0)[lb]{\smash{{\SetFigFont{10}{12.0}{\familydefault}{\mddefault}{\updefault}{\color[rgb]{0,0,0}$h(x_2)$}%
}}}}
\put(5610,-1104){\makebox(0,0)[lb]{\smash{{\SetFigFont{10}{12.0}{\familydefault}{\mddefault}{\updefault}{\color[rgb]{0,0,0}$h(x_4)$}%
}}}}
\put(5097,-950){\makebox(0,0)[b]{\smash{{\SetFigFont{10}{12.0}{\familydefault}{\mddefault}{\updefault}{\color[rgb]{0,0,0}$h(x_3)$}%
}}}}
\put(2529,-280){\makebox(0,0)[b]{\smash{{\SetFigFont{10}{12.0}{\familydefault}{\mddefault}{\updefault}{\color[rgb]{0,0,0}$x_3$}%
}}}}
\put(3095,-1063){\makebox(0,0)[lb]{\smash{{\SetFigFont{10}{12.0}{\familydefault}{\mddefault}{\updefault}{\color[rgb]{0,0,0}$x_4$}%
}}}}
\put(1796,-861){\makebox(0,0)[b]{\smash{{\SetFigFont{10}{12.0}{\familydefault}{\mddefault}{\updefault}{\color[rgb]{0,0,0}$E$}%
}}}}
\put(1878,-474){\makebox(0,0)[b]{\smash{{\SetFigFont{10}{12.0}{\familydefault}{\mddefault}{\updefault}{\color[rgb]{0,0,0}$l$}%
}}}}
\put(4571,-1133){\makebox(0,0)[rb]{\smash{{\SetFigFont{10}{12.0}{\familydefault}{\mddefault}{\updefault}{\color[rgb]{0,0,0}$h(E)$}%
}}}}
\end{picture}%

%% file: pics/ch-complex.tex
\begin{picture}(0,0)%
\includegraphics{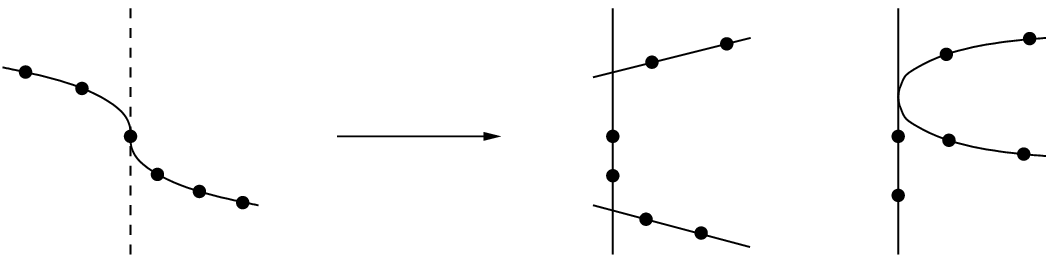}%
\end{picture}%
\setlength{\unitlength}{4144sp}%
\begingroup\makeatletter\ifx\SetFigFont\undefined%
\gdef\SetFigFont#1#2#3#4#5{%
  \reset@font\fontsize{#1}{#2pt}%
  \fontfamily{#3}\fontseries{#4}\fontshape{#5}%
  \selectfont}%
\fi\endgroup%
\begin{picture}(4794,1327)(304,-701)
\put(901,-646){\makebox(0,0)[b]{\smash{{\SetFigFont{10}{12.0}{\familydefault}{\mddefault}{\updefault}{\color[rgb]{0,0,0}$L$}%
}}}}
\put(2206,119){\makebox(0,0)[b]{\smash{{\SetFigFont{10}{12.0}{\familydefault}{\mddefault}{\updefault}{\color[rgb]{0,0,0}specialize}%
}}}}
\put(388,209){\makebox(0,0)[rb]{\smash{{\SetFigFont{10}{12.0}{\familydefault}{\mddefault}{\updefault}{\color[rgb]{0,0,0}$P_6$}%
}}}}
\put(655,128){\makebox(0,0)[rb]{\smash{{\SetFigFont{10}{12.0}{\familydefault}{\mddefault}{\updefault}{\color[rgb]{0,0,0}$P_5$}%
}}}}
\put(868,-67){\makebox(0,0)[rb]{\smash{{\SetFigFont{10}{12.0}{\familydefault}{\mddefault}{\updefault}{\color[rgb]{0,0,0}$P_1$}%
}}}}
\put(1033,-97){\makebox(0,0)[lb]{\smash{{\SetFigFont{10}{12.0}{\familydefault}{\mddefault}{\updefault}{\color[rgb]{0,0,0}$P_2$}%
}}}}
\put(1243,-196){\makebox(0,0)[lb]{\smash{{\SetFigFont{10}{12.0}{\familydefault}{\mddefault}{\updefault}{\color[rgb]{0,0,0}$P_3$}%
}}}}
\put(1447,-253){\makebox(0,0)[lb]{\smash{{\SetFigFont{10}{12.0}{\familydefault}{\mddefault}{\updefault}{\color[rgb]{0,0,0}$P_4$}%
}}}}
\put(3106,-646){\makebox(0,0)[b]{\smash{{\SetFigFont{10}{12.0}{\familydefault}{\mddefault}{\updefault}{\color[rgb]{0,0,0}$L$}%
}}}}
\put(4411,-646){\makebox(0,0)[b]{\smash{{\SetFigFont{10}{12.0}{\familydefault}{\mddefault}{\updefault}{\color[rgb]{0,0,0}$L$}%
}}}}
\put(4006,-16){\makebox(0,0)[b]{\smash{{\SetFigFont{12}{14.4}{\familydefault}{\mddefault}{\updefault}{\color[rgb]{0,0,0}+}%
}}}}
\put(3061,-16){\makebox(0,0)[rb]{\smash{{\SetFigFont{10}{12.0}{\familydefault}{\mddefault}{\updefault}{\color[rgb]{0,0,0}$P_1$}%
}}}}
\put(3061,-196){\makebox(0,0)[rb]{\smash{{\SetFigFont{10}{12.0}{\familydefault}{\mddefault}{\updefault}{\color[rgb]{0,0,0}$P_2$}%
}}}}
\put(4366,-16){\makebox(0,0)[rb]{\smash{{\SetFigFont{10}{12.0}{\familydefault}{\mddefault}{\updefault}{\color[rgb]{0,0,0}$P_1$}%
}}}}
\put(4366,-286){\makebox(0,0)[rb]{\smash{{\SetFigFont{10}{12.0}{\familydefault}{\mddefault}{\updefault}{\color[rgb]{0,0,0}$P_2$}%
}}}}
\end{picture}%

%% file: pics/ch-trop.tex
\begin{picture}(0,0)%
\includegraphics{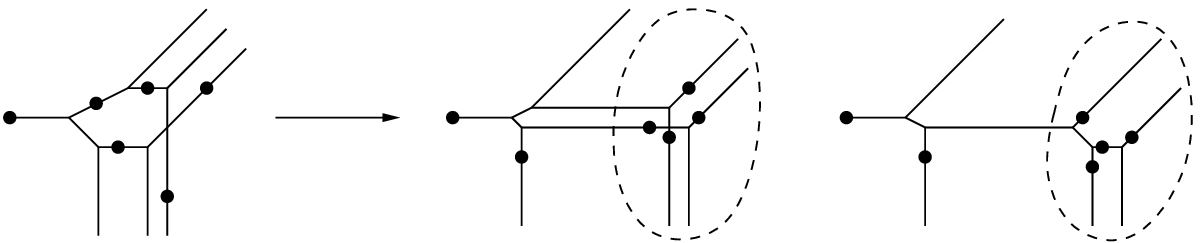}%
\end{picture}%
\setlength{\unitlength}{4144sp}%
\begingroup\makeatletter\ifx\SetFigFont\undefined%
\gdef\SetFigFont#1#2#3#4#5{%
  \reset@font\fontsize{#1}{#2pt}%
  \fontfamily{#3}\fontseries{#4}\fontshape{#5}%
  \selectfont}%
\fi\endgroup%
\begin{picture}(5461,1080)(766,-634)
\put(1321,-331){\makebox(0,0)[b]{\smash{{\SetFigFont{8}{9.6}{\familydefault}{\mddefault}{\updefault}{\color[rgb]{0,0,0}$P_3$}%
}}}}
\put(956,-176){\makebox(0,0)[b]{\smash{{\SetFigFont{8}{9.6}{\familydefault}{\mddefault}{\updefault}{\color[rgb]{0,0,0}$3$}%
}}}}
\put(4366,-106){\makebox(0,0)[b]{\smash{{\SetFigFont{12}{14.4}{\familydefault}{\mddefault}{\updefault}{\color[rgb]{0,0,0}+}%
}}}}
\put(4761,-176){\makebox(0,0)[b]{\smash{{\SetFigFont{8}{9.6}{\familydefault}{\mddefault}{\updefault}{\color[rgb]{0,0,0}$3$}%
}}}}
\put(2966,-176){\makebox(0,0)[b]{\smash{{\SetFigFont{8}{9.6}{\familydefault}{\mddefault}{\updefault}{\color[rgb]{0,0,0}$3$}%
}}}}
\put(5366,-221){\makebox(0,0)[b]{\smash{{\SetFigFont{8}{9.6}{\familydefault}{\mddefault}{\updefault}{\color[rgb]{0,0,0}$2$}%
}}}}
\put(781,-106){\makebox(0,0)[rb]{\smash{{\SetFigFont{8}{9.6}{\familydefault}{\mddefault}{\updefault}{\color[rgb]{0,0,0}$P_1$}%
}}}}
\put(1191, 44){\makebox(0,0)[rb]{\smash{{\SetFigFont{8}{9.6}{\familydefault}{\mddefault}{\updefault}{\color[rgb]{0,0,0}$P_2$}%
}}}}
\put(1446,-66){\makebox(0,0)[rb]{\smash{{\SetFigFont{8}{9.6}{\familydefault}{\mddefault}{\updefault}{\color[rgb]{0,0,0}$P_4$}%
}}}}
\put(1741,-36){\makebox(0,0)[lb]{\smash{{\SetFigFont{8}{9.6}{\familydefault}{\mddefault}{\updefault}{\color[rgb]{0,0,0}$P_6$}%
}}}}
\put(1586,-461){\makebox(0,0)[lb]{\smash{{\SetFigFont{8}{9.6}{\familydefault}{\mddefault}{\updefault}{\color[rgb]{0,0,0}$P_5$}%
}}}}
\put(2801,-106){\makebox(0,0)[rb]{\smash{{\SetFigFont{8}{9.6}{\familydefault}{\mddefault}{\updefault}{\color[rgb]{0,0,0}$P_1$}%
}}}}
\put(3111,-346){\makebox(0,0)[rb]{\smash{{\SetFigFont{8}{9.6}{\familydefault}{\mddefault}{\updefault}{\color[rgb]{0,0,0}$P_2$}%
}}}}
\put(4601,-106){\makebox(0,0)[rb]{\smash{{\SetFigFont{8}{9.6}{\familydefault}{\mddefault}{\updefault}{\color[rgb]{0,0,0}$P_1$}%
}}}}
\put(4961,-341){\makebox(0,0)[rb]{\smash{{\SetFigFont{8}{9.6}{\familydefault}{\mddefault}{\updefault}{\color[rgb]{0,0,0}$P_2$}%
}}}}
\end{picture}%

%% file: summary.tex
\section* {Conclusion}


In the last few years tropical algebraic geometry has evolved with a tremendous
speed. Its general approach to replace algebro-geometric problems by
combinatorial ones often leads to new insights, sometimes to easier proofs of
known statements, and occasionally even to new results in algebraic geometry.

Nevertheless tropical geometry is still in its beginnings since even the most
basic objects of algebraic geometry --- (abstract) varieties and their
morphisms --- do not have a satisfactory counterpart yet in the tropical world.
Consequently, there are many open problems in tropical geometry for the near
future, and one could reasonably expect that the solution to these problems
gives an entirely new strategy to attack many problems in algebraic geometry.
In fact, two such recent new examples in which tropical ideas have already been
applied successfully are the study of compactifications of subvarieties of
algebraic tori (in particular moduli spaces of rational stable curves)
\cite {Te} and low-dimensional topology \cite {Ti}.